\documentclass[a4paper,11pt,twoside,reqno]{amsart}

\usepackage[utf8]{inputenc}
\usepackage[plainpages=false,pdfpagelabels=true]{hyperref}
\usepackage{amssymb,amsthm,wasysym}
\usepackage{slashed}
\usepackage[margin=1in]{geometry}

\newtheorem{Satz}{Theorem}[section]
\newtheorem{Prop}[Satz]{Proposition}
\newtheorem{Lem}[Satz]{Lemma}

\newtheorem{Cor}[Satz]{Corollary}
\newcommand{\cR}{{\mathcal R}}
\theoremstyle{definition}
\newtheorem{Dfn}[Satz]{Definition}
\newtheorem{Bem}[Satz]{Remark}
\newtheorem{Bsp}[Satz]{Example}

\allowdisplaybreaks[1]

\renewcommand{\epsilon}{\varepsilon}
\newcommand{\R}{\ensuremath{\mathbb{R}}}

\newcommand{\D}{\slashed{D}}
\newcommand{\p}{\slashed{\partial}}
\newcommand{\sff}{\mathrm{I\!I}}

\numberwithin{equation}{section}

\title{On the evolution of regularized Dirac-harmonic Maps from closed surfaces}
\author{Volker Branding}
\date{\today}
\address{University of Vienna, Faculty of Mathematics\\
Oskar-Morgenstern-Platz 1, 1090 Vienna, Austria}
\email{volker.branding@univie.ac.at}
\begin{document}

\begin{abstract}
We study the evolution equations for a regularized version of Dirac-harmonic maps from closed Riemannian surfaces.
We establish the existence of a global weak solution for the regularized problem, which is smooth away from finitely many singularities.
Moreover, we discuss the convergence of the evolution equations and address the question if we can remove the regularization in the end.
\end{abstract} 

\maketitle

\section{Introduction and Results}
Harmonic maps from Riemannian surfaces to Riemannian manifolds are a variational problem with rich structure.
Due to their conformal invariance the latter share a lot of special properties.
Among these are for example their regularity and the removal of isolated singularities.
The existence of harmonic maps from surfaces has been established by several methods. 
The approach by Sacks and Uhlenbeck \cite{MR604040} uses a perturbation of the energy functional such that it satisfies the Palais-Smale condition.
The heat flow method was successfully applied in this case by Struwe \cite{MR826871}.

An extension of harmonic maps motivated from supersymmetric field theories in physics are
\emph{Dirac-harmonic maps} introduced in \cite{MR2262709}. These also arise as critical points of an action functional 
and couple the equation for harmonic maps with spinor fields. A Dirac-harmonic map is given by a pair \((\phi,\psi)\) consisting
of a map \(\phi\) and a spinor \(\psi\) along that map. Moreover, Dirac-harmonic maps still belong to the class of conformally invariant variational problems. 
For the physics background see \cite{MR1701618}.

Taking also into account an additional curvature term in the energy functional one is led to \emph{Dirac-harmonic maps with curvature term}, see \cite{MR3333092,MR3558358,MR2370260}.
Dirac-harmonic maps coupled to a two-form potential, called \emph{Magnetic Dirac-harmonic maps}, are studied in \cite{MR3305429}
and Dirac-harmonic maps to manifolds with torsion are examined in \cite{MR3493217}.

At present, many analytical results for Dirac-harmonic maps have already been obtained. 
These include the regularity of solutions \cite{MR2176464,MR2544729,MR2506243}, a removable singularity theorem \cite{MR2262709} and the energy identity \cite{MR2267756}.
Dirac-harmonic maps between closed surfaces are classified in \cite{MR2496649}.
Several vanishing results for Dirac-harmonic maps and their variants can be found in \cite{MR4034775,MR3886921,MR2370260}.

Although many analytical aspects of Dirac-harmonic maps are well understood by now,
the existence question is still not answered in general. 
Some explicit solutions of the Euler-Lagrange equations for Dirac-harmonic maps are given in \cite{MR2569270}.
Using the Atiyah-Singer index theorem uncoupled solutions to the Euler Lagrange equations have been constructed in \cite{MR3070562}.
Namely, for a given map \(\phi_0\) a spinor \(\psi\) is constructed such that the pair \((\phi_0,\psi)\) is a Dirac-harmonic map.
The boundary value problem for Dirac-harmonic maps was treated in \cite{MR3044133} and \cite{MR3085099}.
For a recent survey on mathematical results regarding Dirac-harmonic maps and their variants see \cite{MR3913850}.	

Since Dirac-harmonic maps interpolate between harmonic maps and harmonic spinors,
the existence question for Dirac-harmonic maps can be attacked from two different perspectives. 
On the one hand, one may use methods from spin geometry, as in \cite{MR3070562}, or one may apply methods from the analysis of harmonic maps. 
This of course includes the heat-flow method. However, we cannot apply it directly since the energy functional 
for Dirac-harmonic maps is unbounded from below.

Hence, our approach is to solve an easier problem first and to hope that one can take
a suitable limit in the end. More precisely, we consider the following regularized energy functional
\begin{equation}
\label{energy-regularized}
E_{\epsilon}(\phi,\psi)=\frac{1}{2}\int_M(|d\phi|^2+\langle\psi,\D\psi\rangle+\epsilon|\tilde{\nabla}\psi|^2)dM.
\end{equation}
The first term is the Dirichlet energy of the map \(\phi\), \(\psi\) is a vector spinor and \(\D\) the twisted Dirac operator acting on \(\psi\).
The last term is the \(L^2\)-norm of the covariant derivative of \(\psi\).
Moreover, \(\epsilon>0\) is a bookkeeping parameter.
We study the \(L^2\)-gradient flow of \(E_{\epsilon}(\phi,\psi)\), which is given by
\begin{align}
\label{evolution-phi-surface-intro}\frac{\partial\phi_t}{\partial t}=&\tau(\phi_t)-\cR(\phi_t,\psi_t)-\epsilon\cR_c(\phi_t,\psi_t),\\
\label{evolution-psi-surface-intro}\frac{\tilde{\nabla}\psi_t}{\partial t}=&\epsilon\tilde{\Delta}\psi_t-\D\psi_t
\end{align}
with initial data \((\phi_0,\psi_0)\).
Here, \(\tau(\phi)\) is the tension field of the map \(\phi\) and \(\tilde{\Delta}\) denotes the connection Laplacian for vector spinors.
Since \(\psi\) is a section in the vector bundle \(\Sigma M\otimes\phi^{-1}TN\) we have to use the covariant derivative on this bundle
to calculate the derivative of \(\psi\) with respect to \(t\), which is denoted by \(\frac{\tilde{\nabla}\psi_t}{\partial t}\).
The curvature terms \(\cR(\phi_t,\psi_t)\) and \(\cR_c(\phi_t,\psi_t)\) are of lower order.

Before we state our main result let us make the following observation:
\begin{Bem}
The functional $E_\epsilon(\phi,\psi)$ satisfies
\[
-\frac{1}{4\epsilon}\int_M|\psi|^2dM\leq E_\epsilon(\phi,\psi)\leq\infty.
\]
Thus, the \(L^2\)-norm of \(\psi\) will play an important role in the study of the 
\(L^2\)-gradient flow of the functional \(E_\epsilon(\phi,\psi)\).
\end{Bem}

Our aim is to prove a result similar to Struwe's result \cite{MR826871}, see also \cite{MR965544}, for the harmonic map heat flow from surfaces. 
Due to the coupling between the fields \(\phi\) and \(\psi\) new analytical difficulties arise. Nevertheless, we will prove
\begin{Satz}
\label{theorem-surface}
Let $M$ be a closed Riemannian surface with fixed spin structure and \(N\) a compact Riemannian manifold without boundary.
Suppose that 
\begin{align}
\label{condition-poincare}
\int_M|\psi_t|^2dM\leq c_1\int_M|\tilde\nabla\psi_t|^2dM
\end{align}
holds for all \(t\in [0,\infty)\), where \(c_1>0\).

Then for any smooth initial data $(\phi_0,\psi_0)$ and \(\epsilon>0\) sufficiently large, there exists a global weak solution 
\[
\phi\colon M\times [0,\infty)\to N,\qquad \psi\colon M\times [0,\infty)\to\Sigma M\otimes\phi^{-1}TN
\]
of (\ref{evolution-phi-surface-intro}) and (\ref{evolution-psi-surface-intro}) on
\(M\times [0,\infty)\), which is smooth away from at most finitely many singular points $(x_k,t_k), 1\leq k\leq K$
with \(K=K(\epsilon,\phi_0,\psi_0)\). The weak solution constructed here is unique and the energy functional \eqref{energy-regularized} of the weak solution
is decreasing with respect to time.\\
There exists a sequence \(t_k\to\infty\) such that \((\phi(\cdot,t_k),\psi(\cdot,t_k))\)
converges weakly in $H^{1}(M,N)\times H^{1}(M,\Sigma M\otimes\phi_t^{-1}TN)$ 
to a regularized Dirac-harmonic map \((\phi_\infty,\psi_\infty)\) as $k\to \infty$ suitably
and strongly away from finitely many points \((x_k,t_k=\infty)\).
The pair \((\phi_\infty,\psi_\infty)\) is smooth on \(M\setminus\{x_1,\ldots,x_K\}\).
\end{Satz}

\begin{Bem}
\begin{enumerate}
\item It seems that we have to impose the condition \eqref{condition-poincare} in order to be able to prove Theorem \ref{theorem-surface}.
\item Unfortunately, taking the limit \(\epsilon\to 0\) after \(t\to\infty\) to obtain a Dirac-harmonic map does not seem to be possible.
We will see later, that both the number of singularities and the regularity of \((\phi_\infty,\psi_\infty)\) crucially
depend on \(\epsilon\).
\end{enumerate}

\end{Bem}

A similar approach in the one-dimensional case was performed in \cite{MR3435758},
see also \cite{MR2860404}.
Recently, a new heat-flow approach for Dirac-harmonic maps has been studied in which the Dirac equation
is considered as a constraint while the map is deformed by a heat-type equation.
Several existence results using this approach could be obtained in the
case of a one-dimensional domain \cite{MR3412386} and for the domain being a compact
surface with boundary \cite{MR3724759}.
The short time existence for this flow in the case of a closed manifold was recently established
in \cite{MR3719555}.

The results presented in this article are part of the author’s PhD thesis \cite{phd}.

We would also like to point out that several existence results for Dirac-wave maps could be established \cite{MR3917346,MR3830277,MR2138082}
which are Dirac-harmonic maps from a domain with a Lorentzian metric.

This article is organized as follows. After introducing the framework for Dirac-harmonic maps,
we present a regularized version of Dirac-harmonic maps. Afterwards, we study
the \(L^2\)-gradient flow of the regularized functional in Section 2. 
In Section 3 we establish the existence of a long-time solution and Section 4 then discusses 
the convergence of the evolution equations. In the last section we analyze the limit \(\epsilon\to 0\).

Let us now describe the setup in more detail. We suppose that $M$ is a closed Riemannian spin surface and $N$ a compact Riemannian manifold.
Every orientable Riemannian surface admits a spin structure, the number of different spin structures can be counted by the genus of the surface.
For more details on spin geometry, see the book \cite{MR1031992}. Coordinates on $M$ will be denoted by $x$, whereas coordinates on $N$ will be denoted by $y$.
Indices on $M$ are labeled by Greek letters, whereas indices on $N$ are labeled by Latin letters.
We use the Einstein summation convention, which means that we will sum over repeated indices.

For a given map $\phi:M\to N$, we consider the pull-back bundle $\phi^{-1}TN$ of $TN$
and twist it with the spinor bundle $\Sigma M$. On this twisted bundle $\Sigma M\otimes\phi^{-1}TN$
there is a metric induced from the metrics on $\Sigma M$ and $\phi^{-1}TN$.
The induced connection on $\Sigma M\otimes\phi^{-1}TN$ will be denoted by $\tilde{\nabla}$. 
We will always assume that all connections are metric and free of torsion.
Locally, sections of $\Sigma M\otimes\phi^{-1}TN$, called \emph{vector spinors}, can be expressed as
\[
\psi(x)=\psi^i(x)\otimes\frac{\partial}{\partial y^i}(\phi(x)).
\]
On the spinor bundle \(\Sigma M\) we have the Clifford multiplication of spinors with tangent vectors,
which is skew-symmetric, namely
\[
\langle\psi,X\cdot\chi\rangle_{\Sigma M}=-\langle X\cdot\psi,\chi\rangle_{\Sigma M}
\]
for \(\psi,\chi\in\Gamma(\Sigma M)\) and \(X\in TM\).
We denote the Dirac operator on $\Sigma M$ by $\slashed{\partial}$ and the Dirac operator on the twisted bundle by $\D$,
which is given by
\[
\D=e_\alpha\cdot\tilde{\nabla}_{e_\alpha},
\]
where \(e_\alpha\) is a local basis of \(TM\).
In terms of local coordinates \(\D\psi\) can be expressed as
\[
\D\psi=\slashed{\partial}\psi^i\otimes\frac{\partial}{\partial y^i}(\phi(x))
+\Gamma^i_{jk}\frac{\partial\phi^j}{\partial x_\alpha}e_\alpha\cdot\psi^k(x)\otimes\frac{\partial}{\partial y^i}(\phi(x)),
\]
where \(\Gamma^i_{jk}\) are the Christoffel symbols on \(N\).
Since the connection on $\phi^{-1}TN$ is metric the operator \(\D\) is self-adjoint with respect to the \(L^2\) norm.

We may now state the energy functional for Dirac-harmonic maps
\begin{equation}
\label{energy-dirac-harmonic}
E(\phi,\psi)=\frac{1}{2}\int_M(|d\phi|^2+\langle\psi,\D\psi\rangle)dM,
\end{equation}
which has the critical points (see \cite{MR2262709}, p.\ 413, Prop. 2.1):
\begin{Prop}
The Euler-Lagrange equations for the functional $E(\phi,\psi)$ are given by
\begin{align}
\label{euler-lagrange-phi-unregularized}
\tau(\phi)=&\cR(\phi,\psi), \\
\label{euler-lagrange-psi-unregularized}
\D\psi=&0,
\end{align}
where \(\tau(\phi)\) is the tension field of the map \(\phi\) and the right hand side \(\cR(\phi,\psi)\)
is explicitly given by
\begin{equation}
\label{curvature-term}
\cR(\phi,\psi)=\frac{1}{2}R^N(\psi,e_\alpha\cdot\psi)d\phi(e_\alpha)
\end{equation}
with \(R^N\) being the Riemann curvature tensor on \(N\).
\end{Prop}
In terms of local coordinates, the Euler-Lagrange equations acquire the form
\begin{align*}
\tau^m(\phi)-\frac{1}{2}R^m_{~lij}(\phi)\langle\psi^i,\nabla\phi^l\cdot\psi^j\rangle_{\Sigma M}=&0,\\
\slashed{\partial}\psi^i+\Gamma^i_{jk}(\phi)\frac{\partial\phi^j}{\partial x_\alpha}e_\alpha\cdot\psi^k=&0,
\end{align*}
where \(R^m_{~lij}\) are the components of the curvature tensor on \(N\).
Solutions of the system (\ref{euler-lagrange-phi-unregularized}), (\ref{euler-lagrange-psi-unregularized}) 
are called $\emph{Dirac-harmonic maps}$ from \(M\to N\).

In the analysis of the energy functional $E(\phi,\psi)$ one faces the problem
that it is unbounded from below, since the operator \(\D\) is unbounded. 
To overcome these analytical difficulties, we ``improve'' the energy
functional $E(\phi,\psi)$ by adding a regularizing term, see \eqref{energy-regularized}.
Note that we formally have
\[\lim_{\epsilon\to 0}E_\epsilon(\phi,\psi)=E(\phi,\psi).
\]
Of course, we would like to keep the parameter \(\epsilon\) as small as possible. Unfortunately,
in order to derive energy estimates, we have to drop this assumption.

As a next step we present the Euler-Lagrange equations for \(E_\epsilon(\phi,\psi)\).

\begin{Prop}
The critical points of the functional \(E_\epsilon(\phi,\psi)\) are given by
\begin{align}
\label{euler-lagrange-phi-regularized}
\tau(\phi)=&\cR(\phi,\psi)+\epsilon\cR_c(\phi,\psi), \\
\label{euler-lagrange-psi-regularized}
\epsilon\tilde{\Delta}\psi=&\D\psi
\end{align}
with the curvature term 
\begin{align*}
\cR_c(\phi,\psi)= R^N(\tilde{\nabla}_{e_\alpha}\psi,\psi)d\phi(e_\alpha) \in\Gamma(\phi^{-1}TN)
\end{align*}
and \(\cR(\phi,\psi)\) given by \eqref{curvature-term}.
Moreover, \(\tilde{\Delta}\) denotes the connection Laplacian on the bundle \(\Sigma M\otimes\phi^{-1}TN\).
\end{Prop}

\begin{proof}
For a proof, see \cite{phd}, Section 2.2.
\end{proof}

Written in local coordinates, the new terms arising from the variation of \(E_\epsilon(\phi,\psi)\) acquire the following form:
\begin{align*}
\cR_c(\phi,\psi)=& R^m_{~lij}\frac{\partial}{\partial y^m}\big(\frac{\partial\phi^l}{\partial x_\alpha}\langle\nabla_{e_\alpha}^{\Sigma M}\psi^i,\psi^j\rangle_{\Sigma M} 
+\Gamma^j_{rs}\frac{\partial\phi^l}{\partial x_\alpha}\langle\psi^i,\psi^r\rangle_{\Sigma M}\frac{\partial\phi^s}{\partial x_\alpha}\big),\\
\tilde{\Delta}\psi=&\big(\Delta^{\Sigma M}\psi^m+2\nabla^{\Sigma M}_{e_\alpha}\psi^i\Gamma^m_{ij}\frac{\partial\phi^j}{\partial x_\alpha}
+\psi^i\Gamma^m_{ij,p}\frac{\partial\phi^p}{\partial x_\alpha}\frac{\partial\phi^j}{\partial x_\alpha} 
+\psi^i\Gamma^m_{ij}\frac{\partial^2\phi^j}{\partial x_\alpha^2} \\
&+\psi^i\Gamma_{ij}^k\Gamma^m_{ks}\frac{\partial\phi^j}{\partial x_\alpha}\frac{\partial\phi^s}{\partial x_\alpha}
\big)\otimes\frac{\partial}{\partial y^m}.
\end{align*}
Solutions of the system (\ref{euler-lagrange-phi-regularized}), (\ref{euler-lagrange-psi-regularized}) will be called
\emph{regularized Dirac-harmonic maps} from \(M\to N\).

\begin{Bem}
On a compact Riemann surface the following terms are invariant under conformal transformations:
\[
\int_M|d\phi|^2dM,\qquad \int_M\langle\psi,\D\psi\rangle dM, \qquad \int_M|\psi|^4dM.
\]
A proof can for example be found in \cite{MR2262709}, p.\ 416, Lemma 3.1.
In particular, this means that the functional $E(\phi,\psi)$ is 
conformally invariant in dimension two. 
We will see later that the $L^4$-norm of $\psi$ plays an important role in the context 
of a removable singularity theorem.
On the other hand, we note that through the regularization the conformal invariance is broken.
\end{Bem}

\section{Evolution Equations and Energy Estimates}
We now turn to the \(L^2\)-gradient flow of the regularized functional \(E_\epsilon(\phi,\psi)\):
\begin{align}
\label{evolution-phi-surface}
\frac{\partial\phi_t}{\partial t}=&\tau(\phi_t)-\cR(\phi_t,\psi_t)-\epsilon\cR_c(\phi_t,\psi_t),\\
\label{evolution-psi-surface}
\frac{\tilde{\nabla}\psi_t}{\partial t}=&\epsilon\tilde{\Delta}\psi_t-\D\psi_t
\end{align}
with initial data \((\phi_0,\psi_0)\).

As \(\epsilon>0\) the system \eqref{evolution-phi-surface}, \eqref{evolution-psi-surface} is clearly parabolic
and the existence of a smooth short-time solution up to a time \(T_\textrm{max}\) can be obtained by standard methods,
see Theorem 3.24 in \cite{phd}.

Before turning to the derivation of energy estimates let us make the following remarks.
\begin{Lem}
There does not exist a Dirac-harmonic map from $T^2\to S^2$ with $\deg\phi=\pm 1$.
\end{Lem}
\begin{proof}
The proof is by contradiction. Assume that \((\phi,\psi)\) is a Dirac-harmonic map
from \(T^2\to S^2\) with \(\deg(\phi)=\pm 1\). By the classification theorem for Dirac-harmonic 
maps between surfaces obtained in \cite{MR2496649} the map \(\phi\)
has to be harmonic in this case. On the other hand, Eells and Wood proved in \cite{MR0420708} 
that there does not exist a harmonic map from $T^2\to S^2$ of degree $\pm 1$ independently 
of the metrics chosen on the surfaces $M$ and $N$.
\end{proof}
\begin{Bem}
Since the degree of a map is homotopy-invariant, 
we cannot find a Dirac-harmonic map from \(T^2\to S^2\) in the homotopy class
of \(\phi\) with \(\deg\phi=\pm 1\). This example motivates the occurrence
of singularities in the heat flow for (regularized) Dirac-harmonic maps.
\end{Bem}

\begin{Bem}
We cannot hope to find a global smooth solution of \eqref{evolution-phi-surface} and \eqref{evolution-psi-surface}, as already the harmonic map heat flow
develops singularities in finite time \cite{MR1180392}.
In addition, we cannot expect to find a unique solution in general since in \cite{MR1870959} and \cite{MR1883901}, 
solutions that are different from Struwe's solution \cite{MR826871}, were constructed.
\end{Bem}

In the following we will often need the following combination of quantities
\begin{align*}
E_{\epsilon}(\phi_t,\psi_t,B_R):=&\frac{1}{2}\int_{B_R}(|d\phi_t|^2+\langle\psi_t,\D\psi_t\rangle+\epsilon|\tilde{\nabla}\psi_t|^2)dM,\\
F(\phi_t,\psi_t,B_R):=&\frac{1}{2}\int_{B_R}(|d\phi_t|^2+\epsilon|\tilde{\nabla}\psi_t|^2)dM,\\
F(\phi_t,\psi_t):=&\frac{1}{2}\int_M(|d\phi_t|^2+\epsilon|\tilde{\nabla}\psi_t|^2)dM.
\end{align*}
Moreover, for the further analysis it turns out to be useful to introduce the following function space with \(Q=M\times [0,T)\)
and \(dQ=dMdt\):
\begin{align*}
V:=\bigg\{\sup_{0\leq t\leq T}F(\phi_t,\psi_t)+\int_Q\big(|\nabla^2\phi|^2+|\tilde{\nabla}^2\psi|^2
+\big|\frac{\partial\phi_t}{\partial t}\big|^2+\big|\frac{\tilde{\nabla}\psi_t}{\partial t}\big|^2\big)dQ<\infty 
\bigg\}
\end{align*}

Let \(\Omega\in\R^2\) be a bounded domain.
Then Ladyzhenskaya's inequality holds, that is
\begin{Lem}
Assume that \(v\in H_0^{1}(\Omega)\). Then the following inequality holds:
\begin{equation}
\|v\|^4_{L^4(\Omega)}\leq C\|v\|^2_{L^2(\Omega)}\|\nabla v\|^2_{L^2(\Omega)}
\end{equation}
\end{Lem}

In addition, we need a local version of Ladyzhenskaya's inequality from above.
By \(B_R(x)\) we denote the geodesic ball of radius \(R\) 
around \(x\in M\) and \(i_M\) denotes the injectivity radius of \(M\).
In terms of these quantities we can formulate the following:
\begin{Lem}
Assume that \(v\in V\). Then there exists a constant $C$ such that for any \(R\in(0,i_M)\)
the following inequality holds:
\begin{equation}
\label{local-sobolev-inequality}
\int_M|\nabla v|^4 dM\leq C\sup_{x\in M}\int_{B_{R}(x)\author{}}|\nabla v|^2dM\big(\int_M|\nabla^2v|^2dM+\frac{1}{R^2}\int_M|\nabla v|^2dM\big).
\end{equation}
\end{Lem}
\begin{proof}
A proof can for example be found in \cite{MR2431434}, p.\ 225, Lemma 6.7.
\end{proof}

As a first step, we want to obtain a pointwise bound for the norm of the spinor \(\psi_t\). 
Using (\ref{evolution-psi-surface}) we calculate
\begin{align}
\nonumber \frac{\partial}{\partial t}\frac{1}{2}|\psi_t|^2=&
\frac{\epsilon}{2}\Delta|\psi_t|^2-\langle\psi_t,\D\psi_t\rangle-\epsilon|\tilde{\nabla}\psi_t|^2 \\
\nonumber \leq& \frac{\epsilon}{2}\Delta|\psi_t|^2+\sqrt{2}|\psi_t||\tilde{\nabla}\psi_t|-\epsilon|\tilde{\nabla}\psi_t|^2 \\
\leq& \frac{\epsilon}{2}\Delta|\psi_t|^2+\frac{1}{2\epsilon}|\psi_t|^2.
\label{evolution-pointwise-spinor}
\end{align}

\begin{Bem}
If we apply the maximum principle to \eqref{evolution-pointwise-spinor}, we obtain the estimate
\[
|\psi_t|^2\leq|\psi_0|^2e^{\frac{t}{\epsilon}}.
\]
In particular, if the initial spinor \(\psi_0\) vanishes, then
our system \eqref{evolution-phi-surface} and \eqref{evolution-psi-surface} reduces to
the harmonic map heat flow studied by Struwe in \cite{MR826871}. Moreover, if \(\psi_t=0\)
for some time \(t\) then \(\psi_t=0\) for all \(T\geq t\).
\end{Bem}

\begin{Lem}
\label{lemma-bound-psi}
Let \(\psi_t\in C^2(M\times[0,T),\Sigma M\otimes\phi_t^{-1}TN)\) be a solution of (\ref{evolution-psi-surface})
and assume that that \eqref{condition-poincare} holds.
For \(\epsilon\) large enough we get a uniform bound on \(\psi_t\) 
\begin{equation}
|\psi_t|^2_{L^\infty(M\times[0,T))}\leq Ce^\frac{1}{\epsilon}.
\end{equation}
The constant \(C\) depends on \(M,N,c_1\) and the \(L^2\)-norm of \(\psi_0\).
\end{Lem}
\begin{proof}
We already know that \(\psi_t\) solves the pointwise equation \eqref{evolution-pointwise-spinor}.
If we can also bound the \(L^2\)-norm of the spinor \(\psi_t\) we get a uniform pointwise bound 
by Lemma \ref{maximum-principle-l2} (see the appendix for its precise formulation). Thus, we calculate
\begin{align*}
\frac{\partial}{\partial t}\frac{1}{2}\int_M|\psi_t|^2dM&=
-\int_M\langle\psi_t,\D\psi_t\rangle dM-\epsilon\int_M|\tilde{\nabla}\psi_t|^2dM \\
&\leq -\frac{\epsilon}{2}\int|\tilde{\nabla}\psi_t|^2dM+\frac{1}{2\epsilon}\int_M|\psi_t|^2dM\\
&\leq \big(-\frac{\epsilon}{2}+\frac{c_1}{2\epsilon}\big)\int|\tilde\nabla\psi_t|^2dM,
\end{align*}
where we applied \eqref{condition-poincare} in the last step. Hence for \(\epsilon>0\) big enough
the right hand side of the above equation will be negative which gives the desired bound
on the \(L^2\)-norm of \(\psi_t\).
\end{proof}

In the following \(C\) denotes a universal constant that may change from line to line.
Since our evolution equations are originating from a variational problem, 
we get bounds in terms of the initial data \((\phi_0,\psi_0)\).
\begin{Lem}
\label{lemma-surface-estimates-initial-data}
Let \((\phi_t,\psi_t)\in V\) be a solution of (\ref{evolution-phi-surface}) and (\ref{evolution-psi-surface}).
If in addition \(\int_M|\psi_t|^2dM\leq C\), then we have for all \(t\in [0,T)\)
\[
\int_M(|d\phi_t|^2+\epsilon|\tilde{\nabla}\psi_t|^2)dM+
\int_Q \big(\big|\frac{\partial\phi_t}{\partial t}\big|^2+\big|\frac{\tilde{\nabla}\psi_t}{\partial t}\big|^2\big)dMdt\leq C.
\]
The constant $C$ depends on \(M,\epsilon,E_\epsilon(\phi_0,\psi_0)\) and \(\psi_0\).
\end{Lem}

\begin{proof}
The inequality follows from the fact that the system (\ref{evolution-phi-surface}), (\ref{evolution-psi-surface}) is the 
\(L^2\)-gradient flow of the functional \(E_\epsilon(\phi,\psi)\) and
\[
-\frac{1}{2}\langle\psi_t,\D\psi_t\rangle\leq \frac{\epsilon}{8}|\tilde{\nabla}\psi_t|^2+\frac{1}{2\epsilon}|\psi_t|^2.
\]
\end{proof}

The next Lemma is the analogue of Lemma 3.6 from \cite{MR826871}.
We want to get local bounds of the $L^2$-norms of $d\phi_t$ and $\tilde{\nabla}\psi_t$.
\begin{Lem}
\label{lemma-energy-local-surface}
Let \((\phi_t,\psi_t)\in V\) be a solution of (\ref{evolution-phi-surface}) and (\ref{evolution-psi-surface}).
For $R\in (0,i_M)$ and any $(x,t)\in Q$ there holds the estimate
\begin{equation}
E_{\epsilon}(\phi_t,\psi_t,B_R)\leq\frac{C}{R^2}\int_Q(|d\phi_t|^2+|\psi_t|^2+\epsilon^2|\tilde{\nabla}\psi_t|^2)dQ+E_{\epsilon}(\phi_0,\psi_0,B_{2R}),
\end{equation}
where the constant $C$ only depends on $M$.
\end{Lem}
\begin{proof}
First of all, we choose a smooth cut-off function $\eta$ with the following properties
\begin{eqnarray*}
\eta\in C^\infty(M),\qquad\eta\geq 0,\qquad\eta=1~\textrm{on}~B_R(x_0),\\
\eta=0~\textrm{on}~M\setminus B_{2R}(x_0),\qquad |\nabla\eta|_{L^\infty}\leq\frac{C}{R},
\end{eqnarray*}
where again \(B_R(x_0)\) denotes the geodesic ball of radius \(R\) around \(x_0\in M\).
In addition, we choose an orthonormal basis \(\{e_\alpha,\alpha=1,2\}\) on \(M\)
such that \(\nabla_{e_\alpha}e_\beta=\nabla_{\partial_t}e_\alpha=0\) at the considered point.
By a direct calculation we find
\begin{align*}
\frac{\partial}{\partial t}\frac{1}{2}|d\phi_t|^2=&\partial_{e_\alpha}\langle\frac{\partial\phi_t}{\partial t},d\phi_t(e_\alpha)\rangle
-\langle\frac{\partial\phi_t}{\partial t},\tau(\phi_t)\rangle,\\
\frac{\partial}{\partial t}\frac{1}{2}\langle\psi_t,\D\psi_t\rangle
=&-\partial_{e_\alpha}\frac{1}{2}\langle\frac{\tilde{\nabla}\psi_t}{\partial t},e_\alpha\cdot\psi_t\rangle
+\langle\frac{\tilde{\nabla}\psi_t}{\partial t},\D\psi_t\rangle+\langle\frac{\partial\phi_t}{\partial t},\cR(\phi_t,\psi_t)\rangle,\\ 
\frac{\partial}{\partial t}\frac{1}{2}|\tilde{\nabla}\psi_t|^2=&\langle \cR_c(\phi_t,\psi_t),\frac{\partial\phi_t}{\partial t}\rangle
+\partial_{e_\alpha}\langle\frac{\tilde{\nabla}\psi_t}{\partial t},\tilde{\nabla}_{e_\alpha}\psi_t\rangle
-\langle\frac{\tilde{\nabla}\psi_t}{\partial t},\tilde{\nabla}_{e_\alpha}\tilde{\nabla}_{e_\alpha}\psi_t\rangle.
\end{align*}
Multiplying each of the terms with the cut-off function $\eta^2$, adding up the three terms and
using the evolution equations (\ref{evolution-phi-surface}) and (\ref{evolution-psi-surface}), we find
\begin{align*}
\frac{\partial}{\partial t}\frac{1}{2}\int_M\eta^2(|d\phi_t|^2+\langle\psi_t,\D\psi_t\rangle+\epsilon|\tilde{\nabla}\psi_t|^2)dM
+\int_M\eta^2\big(\big|\frac{\tilde{\nabla}\psi_t}{\partial t}\big|^2
+\big|\frac{\partial\phi_t}{\partial t}\big|^2\big)dM\\
=\int_M\eta^2\partial_{e_\alpha}\big(
\langle\frac{\partial\phi_t}{\partial t},d\phi_t(e_\alpha)\rangle
-\frac{1}{2}\langle\frac{\tilde{\nabla}\psi_t}{\partial t},e_\alpha\cdot\psi_t\rangle
-\epsilon\langle\frac{\tilde{\nabla}\psi_t}{\partial t},\tilde{\nabla}_{e_\alpha}\psi_t\rangle
\big)dM.
\end{align*}
Using integration by parts we derive
\begin{align*}
\int_M\eta^2\partial_{e_\alpha}\langle\frac{\partial\phi_t}{\partial t},d\phi_t(e_\alpha)\rangle dM
\leq&C\int_M|\eta||\nabla\eta||\frac{\partial\phi_t}{\partial t}||d\phi_t|dM, \\
\int_M\eta^2\partial_{e_\alpha}\langle\frac{\tilde{\nabla}\psi_t}{\partial t},e_\alpha\cdot\psi_t\rangle dM
\leq&C\int_M|\eta||\nabla\eta||\frac{\tilde{\nabla}\psi_t}{\partial t}||\psi_t|dM,\\
\int_M\eta^2\partial_{e_\alpha}\langle\frac{\tilde{\nabla}\psi_t}{\partial t},\tilde{\nabla}_{e_\alpha}\psi_t\rangle dM
\leq&C\int_M|\eta||\nabla\eta||\frac{\tilde{\nabla}\psi_t}{\partial t}||\tilde{\nabla}\psi_t|dM.
\end{align*}
Applying Young's inequality and by the properties of the cut-off function $\eta$, we find
\[
\frac{\partial}{\partial t} E_{\epsilon}(\phi_t,\psi_t,B_R)\leq\frac{C}{R^2}\int_M(|d\phi_t|^2+|\psi_t|^2+\epsilon^2|\tilde{\nabla}\psi_t|^2)dM.
\]
Integration with respect to $t$ yields the result.
\end{proof}
We can use the previous Lemma to formulate monotonicity formulas for \(F(\phi_t,\psi_t,B_R)\).
By Young's inequality and the ``monotonicity formula'' for the local energy \(E_\epsilon(\phi_t,\psi_t,B_R)\), we get
\begin{equation}
\label{surface-inequality-F-local}
F(\phi_t,\psi_t,B_R)\leq 2E_\epsilon(\phi_0,\psi_0,B_{2R})+C\frac{T}{R^2}
+\frac{1}{2\epsilon}\int_{B_R}|\psi_t|^2 dM.
\end{equation}
Roughly speaking, we want to make the left hand side of this inequality as small as we have to.
This can be achieved by choosing the initial data \((\phi_0,\psi_0)\), the radius \(R\)
of the ball \(B_R\) and the time \(T\) appropriately. More precisely, we get the following
\begin{Cor}
Let \((\phi_t,\psi_t)\in V\) be a solution of (\ref{evolution-phi-surface}) and (\ref{evolution-psi-surface}).
For a positive constant $\delta_1$ there exist $R\in(0,i_M)$ and \(T_1>0\)
with \(|\psi_t|_{L^\infty(M\times [0,T_1))}\leq C\) such that
\begin{equation}
\sup_{\genfrac{}{}{0pt}{}{x\in M}{0\leq t\leq T_1}} F(\phi_t,\psi_t,B_R)<\delta_1.
\end{equation}
\end{Cor}
\begin{proof}
From Lemma \ref{lemma-energy-local-surface} and the bound on the norm of \(\psi_t\), it follows 
that for any $\delta_1$ and $(\phi_0,\psi_0)$ suitably, there exists a number $R>0$
for which
\begin{equation}
\label{F-local-proof-1}
\sup_{x\in M} \big(2E_\epsilon(\phi_0,\psi_0,B_{2R})+\frac{1}{2\epsilon}\int_{B_R}|\psi_t|^2dM\big)<\frac{\delta_1}{2}.
\end{equation}
For $T_1=\frac{\delta_1R^2}{2C}$ we then get
\begin{equation}
\label{F-local-proof-2}
\sup_{\genfrac{}{}{0pt}{}{x\in M}{0\leq t\leq T_1}} F(\phi_t,\psi_t,B_R(x))<\delta_1,
\end{equation}
such that the desired estimate holds.
\end{proof}

In order to turn the Laplace type terms into full second derivatives, we will
make use of the following Bochner type formulas:
\begin{Lem}[Bochner type formulas]
For a map \(\phi:M\to N\) and a vector spinor \(\psi\in\Gamma(\Sigma M\otimes\phi^{-1}TN)\) 
the following Bochner type formulas hold:
\begin{align}
\label{lemma-surface-bochner-phi} \int_M|\tau(\phi)|^2dM=&\int_M(|\nabla d\phi|^2
+\langle d\phi(Ric^M(e_\beta)),d\phi(e_\beta)\rangle \\
\nonumber &-\langle R^N(d\phi(e_\alpha),d\phi(e_\beta))d\phi(e_\beta),d\phi(e_\alpha)\rangle) dM, \\
\label{lemma-surface-bochner-psi} \int_M|\tilde{\Delta}\psi|^2dM=&\int_M(|\tilde{\nabla}^2\psi|^2
+\langle R^{E_1}(e_\alpha,e_\beta)\tilde{\nabla}_{e_\beta}\psi,\tilde{\nabla}_{e_\alpha}\psi\rangle_{E_1} 
+\langle R^{E_2}(e_\alpha,e_\beta)\psi, \tilde{\nabla}_{e_\beta}\tilde{\nabla}_{e_\alpha}\psi\rangle_{E_2})dM
\end{align}
with the vector bundles \(E_1=T^*M\otimes\Sigma M\otimes\phi^{-1}TN\) and \(E_2=T^*M\otimes E_1\).
\end{Lem}
\begin{proof}
This follows from a direct calculation.
\end{proof}

We are now able to bound the \(L^2\)-norm of the second derivatives of \((\phi_t,\psi_t)\) on \(M\times[0,T_1)\).
\begin{Prop}
Let \((\phi_t,\psi_t)\in V\) be a solution of (\ref{evolution-phi-surface}) and (\ref{evolution-psi-surface}).
Choose \(R>0\) and \(T_1>0\) such that \eqref{F-local-proof-1} and \eqref{F-local-proof-2} hold.
Moreover, assume that \(|\psi_t|_{L^\infty(M\times [0,T_1))}\leq C\). 
Then we have for all $t\in[0,T_1)$
\begin{equation}
\int_Q\big(|\nabla d\phi_t|^2+\epsilon^2|\tilde{\nabla}^2\psi_t|^2\big)dMdt\leq C\big(1+\frac{T_1}{R^2}\big),
\end{equation}
where the constant $C$ depends on $M,N,\epsilon,\psi_0,d\phi_0\) and \(\tilde{\nabla}\psi_0$.
\end{Prop}

\begin{proof}
Using the evolution equations (\ref{evolution-phi-surface}) and (\ref{evolution-psi-surface}) we compute
\begin{align*}
\frac{\partial}{\partial t}\frac{1}{2}\int_M(&|d\phi_t|^2+\epsilon|\tilde{\nabla}\psi_t|^2)dM
+\int_M(|\tau(\phi_t)|^2+\epsilon^2|\tilde{\Delta}\psi_t|^2)dM\\
=&\int_M(\epsilon\langle\cR(\phi_t,\psi_t),\cR_c(\phi_t,\psi_t)\rangle
+\epsilon\langle\tilde{\Delta}\psi_t,\D\psi_t\rangle-\epsilon^2|\cR_c(\phi_t,\psi_t)|^2)dM\\
&+\int_M(\langle\tau(\phi_t),\cR(\phi_t,\psi_t)+2\epsilon\cR_c(\phi_t,\psi_t)\rangle)dM.
\end{align*}
Applying Young's inequality and estimating the terms on the right hand side, we get
\begin{align*}
\frac{\partial}{\partial t}\frac{1}{2}\int_M(&|d\phi_t|^2+\epsilon|\tilde{\nabla}\psi_t|^2)dM
+\frac{1}{2}\int_M|\tau(\phi_t)|^2+\epsilon^2|\tilde{\Delta}\psi_t|^2)dM \\
&\leq\int_M\big(\frac{1}{2}|\cR(\phi_t,\psi_t)|^2+\epsilon\langle\cR(\phi_t,\psi_t),\cR_c(\phi_t,\psi_t)\rangle
+\epsilon^2|\cR_c(\phi_t,\psi_t)|^2 
+\frac{1}{2}|\D\psi_t|^2\big)dM \\
&\leq C\int_M(|d\phi_t|^2|\psi_t|^4+\epsilon^2|d\phi_t|^2|\tilde{\nabla}\psi_t|^2|\psi_t|^2+|\tilde{\nabla}\psi_t|^2)dM\\
&\leq C\int_M(|d\phi_t|^2+\epsilon^2|d\phi_t|^2|\tilde{\nabla}\psi_t|^2+|\tilde{\nabla}\psi_t|^2)dM.
\end{align*}
As a next step we transform the Laplace type terms into second derivatives,
therefore we apply the Bochner type formulas (\ref{lemma-surface-bochner-phi}), (\ref{lemma-surface-bochner-psi}) and find
\begin{align*}
\frac{\partial}{\partial t}\frac{1}{2}\int_M(&|d\phi_t|^2+\epsilon|\tilde{\nabla}\psi_t|^2)dM
+C\int_M(|\nabla d\phi_t|^2+\epsilon^2|\tilde{\nabla}^2\psi_t|^2)dM \\
&\leq C\big(\int_M(|d\phi_t|^2+|\tilde{\nabla}\psi_t|^2)dM
+\int_M(|d\phi_t|^4+\epsilon^2|\tilde{\nabla}\psi_t|^4)dM\big),
\end{align*}
where we estimated all curvature contributions.
Finally, we apply the local Sobolev inequality (\ref{local-sobolev-inequality}) to $\int_M|d\phi_t|^4dM$ and $\int_M|\tilde{\nabla}\psi_t|^4dM$,
which leads to
\begin{align*}
\frac{\partial}{\partial t}\frac{1}{2}\int_M(&|d\phi_t|^2+\epsilon|\tilde{\nabla}\psi_t|^2)dM
+C\int_M(|\nabla d\phi_t|^2+\epsilon^2|\tilde{\nabla}^2\psi_t|^2)dM \\
\leq & C\big(\int_M(|d\phi_t|^2+|\tilde{\nabla}\psi_t|^2)dM
+\frac{\delta_1}{R^2}\int_M(|d\phi_t|^2+\epsilon|\tilde{\nabla}\psi_t|^2)dM\\
&+\delta_1\int_M(|\nabla d\phi_t|^2+\epsilon^2|\tilde{\nabla}^2\psi_t|^2)dM\big).
\end{align*}
Choosing $\delta_1$ small enough, the terms containing the second derivatives on the
right hand side can be absorbed into the left hand side.
Integrating with respect to $t$ yields the result.
\end{proof}

Using the bounds on the second derivatives, we can apply the Sobolev embedding theorem
to bound $\int_Q|d\phi_t|^4dQ$ and $\int_Q|\tilde{\nabla}\psi_t|^4dQ$.
\begin{Cor}
Let \((\phi_t,\psi_t)\in V\) be a solution of (\ref{evolution-phi-surface}) and (\ref{evolution-psi-surface}) with
\(|\psi_t|_{L^\infty(M\times [0,T_1))}\leq C\).  If \(\sup_{(x,t)\in M\times[0,T_1)}F(\phi_t,\psi_t,B_R(x))<\delta_1\), 
then we have for all $t\in[0,T_1)$
\begin{align}
\int_Q|d\phi_t|^4dQ\leq Cf_1(t),\qquad
\int_Q|\tilde{\nabla}\psi_t|^4dQ\leq Cf_2(t)
\end{align}
with \(f_i(t)\) satisfying \(f_i(t)\to 0\) as \(t\to 0\) for \(i=1,2\).
\end{Cor}
\begin{proof}
The bounds follow from the Sobolev embedding in two dimensions and the previous estimates,
namely
\begin{align*}
\int_Q|d\phi_t|^4 dQ\leq C\int_Q|d\phi_t|^2dQ\int_Q|\nabla d\phi_t|^2dQ
\leq Cf_1(t). 
\end{align*}
The estimate on \(\int_Q|\tilde{\nabla}\psi_t|^4dQ\) can be derived by the same method.
\end{proof}

\begin{Cor}
\label{surface-corollay-small-l4}
For \(\delta_1\) small enough and integrating over a small time interval \(|t-s|\leq\delta_2\),
we can achieve
\begin{equation}
\int_s^t\int_M|d\phi_t|^4dQ \leq C,\qquad \int_s^t\int_M|\tilde{\nabla}\psi_t|^4dQ\leq C
\end{equation}
and the right hand side can be made as small as needed. \\
The constant $C$ depends on $M,N,R,\delta_1,\delta_2,\epsilon,\psi_0,d\phi_0\) and \(\tilde{\nabla}\psi_0$.
\end{Cor}

So far, we have derived integral estimates on \(Q=M\times [0,T_1)\) of the second derivatives.
In order to turn these into estimates on \(M\), we have to gain control over the derivatives
with respect to \(t\) of the pair \((\phi_t,\psi_t)\). 
\begin{Lem}
\label{lemma-surface-evolution-kinetic-energies}
Let \((\phi_t,\psi_t)\in V\) be a solution of (\ref{evolution-phi-surface}) and (\ref{evolution-psi-surface}).
Then we have 
\begin{align}
\frac{\partial}{\partial t}\frac{1}{2}\int_M\big|\frac{\partial\phi_t}{\partial t}\big|^2dM=
\int_M&\big(-|\nabla\frac{\partial\phi_t}{\partial t}|^2
+\langle R^N(d\phi_t(e_\alpha),\frac{\partial\phi_t}{\partial t})\frac{\partial\phi_t}{\partial t},d\phi_t(e_\alpha)\rangle\\ \nonumber 
&-\langle\frac{\nabla}{\partial t} \cR(\psi_t,\phi_t),\frac{\partial\phi_t}{\partial t}\rangle
-\epsilon\langle\frac{\nabla}{\partial t}\cR_c(\psi_t,\phi_t),\frac{\partial\phi_t}{\partial t}\rangle\big)dM,\\
\frac{\partial}{\partial t}\frac{1}{2}\int_M\big|\frac{\tilde{\nabla}\psi_t}{\partial t}\big|^2dM
=\int_M&\big(-\epsilon\big|\tilde{\nabla}\frac{\tilde{\nabla}\psi_t}{\partial t}\big|^2
  +\epsilon\langle\frac{\tilde{\nabla}\psi_t}{\partial t},R^N(\frac{\partial\phi_t}{\partial t},d\phi_t(e_\alpha))\tilde{\nabla}_{e_\alpha}\psi_t\rangle\\
\nonumber & -\langle\frac{\tilde{\nabla}\psi_t}{\partial t},\frac{\tilde{\nabla}}{\partial t}\D\psi_t\rangle 
-\epsilon\langle\tilde{\nabla}_{e_\alpha}\frac{\tilde{\nabla}\psi_t}{\partial t},R^{N}(\frac{\partial\phi_t}{\partial t},d\phi_t(e_\alpha))\psi_t\rangle\big)dM
\end{align}
for all \(t\in[0,T)\).
\end{Lem}

\begin{proof}
This follows by a direct calculation.
\end{proof}

\begin{Prop}
\label{lemma-surface-bounding-kinetic-energies}
Let \((\phi_t,\psi_t)\in V\) be a solution of (\ref{evolution-phi-surface}) and (\ref{evolution-psi-surface}) with
\(|\psi_t|_{L^\infty(M\times [0,T_1))}\leq C\).
If 
\(\sup_{(x,t)\in M\times[0,T_1)}F(\phi_t,\psi_t,B_R(x))<\delta_1\) is small enough, we find for \(\tau>0\)
\begin{equation}
\sup_{2\tau\leq t\leq T_1}\int_M\big(|\frac{\partial\phi(\cdot,t)}{\partial t}|^2+|\frac{\tilde{\nabla}\psi(\cdot,t)}{\partial t}|^2\big)dM
\leq C(1+\tau^{-1}),
\end{equation}
where the constant \(C\) depends on $M,N,R,\delta_1,\delta_2,\epsilon,\tau,\psi_0,d\phi_0\) and \(\tilde{\nabla}\psi_0$.
\end{Prop}

\begin{proof}
First of all, we choose an orthonormal basis \(\{e_\alpha,\alpha=1,2\}\) on \(M\)
such that \(\nabla_{\partial_t}e_\alpha=0\) at a considered point.
Combining both equations from Lemma \ref{lemma-surface-evolution-kinetic-energies}
we get
\begin{align*}
\frac{\partial}{\partial t}\frac{1}{2}\int_M\big(\big|&\frac{\partial\phi_t}{\partial t}\big|^2
+\big|\frac{\tilde{\nabla}\psi_t}{\partial t}\big|^2\big)dM+\int_M\big(|\nabla\frac{\partial\phi_t}{\partial t}|^2
+\epsilon\big|\tilde{\nabla}\frac{\tilde{\nabla}\psi_t}{\partial t}\big|^2\big)dM\\&=
\int_M\big(\langle R^N(d\phi_t(e_\alpha),\frac{\partial\phi_t}{\partial t})\frac{\partial\phi_t}{\partial t},d\phi_t(e_\alpha)\rangle
-\langle\frac{\nabla}{\partial t} \cR(\psi_t,\phi_t),\frac{\partial\phi_t}{\partial t}\rangle\\
&-\epsilon\langle\frac{\nabla}{\partial t}\cR_c(\psi_t,\phi_t),\frac{\partial\phi_t}{\partial t}\rangle
+\epsilon\langle\frac{\tilde{\nabla}\psi_t}{\partial t},R^N(\frac{\partial\phi_t}{\partial t},d\phi_t(e_\alpha))\tilde{\nabla}_{e_\alpha}\psi_t\rangle\\
&-\epsilon\langle\tilde{\nabla}_{e_\alpha}\frac{\tilde{\nabla}\psi_t}{\partial t},R^{N}(\frac{\partial\phi_t}{\partial t},d\phi_t(e_\alpha))\psi_t\rangle
-\langle\frac{\tilde{\nabla}\psi_t}{\partial t},\frac{\tilde{\nabla}}{\partial t}\D\psi_t\rangle\big)dM\\
&=A_1+A_2+A_3+A_4+A_5+A_6.
\end{align*}
We have to estimate all terms on the right hand side, starting with the \(A_1\) term
\[
\langle R^N(d\phi_t(e_\alpha),\frac{\partial\phi_t}{\partial t})\frac{\partial\phi_t}{\partial t},d\phi_t(e_\alpha)\rangle 
\leq C|d\phi_t|^2|\frac{\partial\phi_t}{\partial t}|^2.
\]
Calculating directly using the fact that $\psi_t$ is bounded uniformly we find for the \(A_2\) term
\begin{align*}
|\langle\frac{\nabla}{\partial t}\cR(\phi_t,\psi_t),\frac{\partial\phi_t}{\partial t}\rangle|
&\leq C\big(|\frac{\partial\phi_t}{\partial t}|^2|d\phi_t||\psi_t|^2+|d\phi_t||\frac{\partial\phi_t}{\partial t}||\frac{\tilde{\nabla}\psi_t}{\partial t}||\psi_t| 
+|\frac{\nabla}{\partial t}d\phi_t||\frac{\partial\phi_t}{\partial t}||\psi_t|^2\big) \\
&\leq C\big(|d\phi_t|^2|\frac{\partial\phi_t}{\partial t}|^2+|\frac{\partial\phi_t}{\partial t}|^2
+|\frac{\tilde{\nabla}\psi_t}{\partial t}|^2|d\phi_t|^2\big)+\frac{1}{8}|\frac{\nabla}{\partial t}d\phi_t|^2.
\end{align*}

Performing the same manipulations with the \(A_3\) term, we get
\begin{align*}
|\langle\frac{\nabla}{\partial t}\cR_c(\phi_t,\psi_t),\frac{\partial\phi_t}{\partial t}\rangle|
\leq& C\big(|\frac{\partial\phi_t}{\partial t}|^2|d\phi_t||\tilde{\nabla}\psi_t||\psi_t|+
|\tilde{\nabla}\psi_t||d\phi_t||\frac{\partial\phi_t}{\partial t}||\frac{\tilde{\nabla}\psi_t}{\partial t}|\\
&+|\frac{\tilde{\nabla}}{\partial t}\tilde{\nabla}\psi_t||d\phi_t||\frac{\partial\phi_t}{\partial t}||\psi_t| 
+|\tilde{\nabla}\psi_t||\frac{\nabla}{\partial t}d\phi_t||\frac{\partial\phi_t}{\partial t}||\psi_t|\big)\\
\leq &C(|d\phi_t|^2|\frac{\partial\phi_t}{\partial t}|^2+|\tilde{\nabla}\psi_t|^2|\frac{\partial\phi_t}{\partial t}|^2
+|\frac{\tilde{\nabla}\psi_t}{\partial t}|^2|d\phi_t|^2) \\
&+\frac{1}{8}|\frac{\nabla}{\partial t}d\phi_t|^2+\frac{\epsilon}{8}|\frac{\tilde{\nabla}}{\partial t}\tilde{\nabla}\psi_t|^2.
\end{align*}

As a next step, we want to control the terms arising from interchanging covariant spinorial derivatives,
namely \(A_4,A_5\) and \(A_6\).
\begin{align*}
A_4&\leq C|d\phi_t||\frac{\partial\phi_t}{\partial t}||\frac{\tilde{\nabla}\psi_t}{\partial t}||\tilde{\nabla}\psi_t|
\leq C(|d\phi_t|^2|\frac{\partial\phi_t}{\partial t}|^2+|\frac{\tilde{\nabla}\psi_t}{\partial t}|^2|\tilde{\nabla}\psi_t|^2),\\
A_5&\leq C|d\phi_t||\frac{\partial\phi_t}{\partial t}||\tilde{\nabla}\frac{\tilde{\nabla}\psi_t}{\partial t}|
\leq C|d\phi_t|^2|\frac{\partial\phi_t}{\partial t}|^2+\frac{\epsilon}{8}|\tilde{\nabla}\frac{\tilde{\nabla}\psi_t}{\partial t}|^2.
\end{align*}
Regarding \(A_6\), we use the pointwise bound on \(\psi_t\), interchange covariant derivatives,
estimate the curvature terms and find
\begin{align*}
|\langle\frac{\tilde{\nabla}\psi_t}{\partial t},&\frac{\tilde{\nabla}}{\partial t}\D\psi_t\rangle| \leq
C|\frac{\partial\phi_t}{\partial t}|^2|d\phi_t|^2+C|\frac{\tilde{\nabla}\psi_t}{\partial t}|^2
+\frac{\epsilon}{8}|\tilde{\nabla}\frac{\tilde{\nabla}\psi_t}{\partial t}|^ 2.
\end{align*}
Note that \(\nabla\frac{\partial\phi}{\partial t}=\frac{\nabla}{\partial t}d\phi\), which
is due to the torsion freeness of the connection.
We sum up the different contributions and find the following inequality
\begin{align*}
\frac{\partial}{\partial t}\frac{1}{2}&\int_M\big(\big|\frac{\partial\phi_t}{\partial t}\big|^2
+\big|\frac{\tilde{\nabla}\psi_t}{\partial t}\big|^2\big)dM
+\frac{1}{2}\int_M\big( \big|\nabla\frac{\partial\phi_t}{\partial t}\big|^2
+\epsilon\big|\tilde{\nabla}\frac{\tilde{\nabla}\psi_t}{\partial t}\big|^2\big)dM  \\
\leq &C\big(\int_M|d\phi_t|^2|\frac{\partial\phi_t}{\partial t}|^2dM
+\int_M|\frac{\partial\phi_t}{\partial t}|^2|\tilde{\nabla}\psi_t|^2dM\\
&\hspace{0.4cm}+\int_M|\frac{\tilde{\nabla}\psi_t}{\partial t}|^2|\tilde{\nabla}\psi_t|^2dM
+\int_M\big(|\frac{\partial\phi_t}{\partial t}|^2+|\frac{\tilde{\nabla}\psi_t}{\partial t}|^2\big)dM\big).
\end{align*}
We used part of the second order terms on the left hand side to absorb the second order terms
from the right hand side.\\
Integrating with respect to \(t\) over the domain \(\tau\leq s<t\leq T\) we get
\begin{align*}
\int_s^tdt\frac{\partial}{\partial t}&\frac{1}{2}\int_M\big(\big|\frac{\partial\phi_t}{\partial t}\big|^2+\big|\frac{\tilde{\nabla}\psi_t}{\partial t}|^2\big)dM
+\frac{1}{2}\int_s^t\int_M\big(\big|\nabla\frac{\partial\phi_t}{\partial t}\big|^2
+\epsilon\big|\tilde{\nabla}\frac{\tilde{\nabla}\psi_t}{\partial t}\big|^2\big)dQ \\
\leq & 
C\big(\int_s^t\int_M|d\phi_t|^2|\frac{\partial\phi_t}{\partial t}|^2dQ
+\int_s^t\int_M|\frac{\partial\phi_t}{\partial t}|^2|\tilde{\nabla}\psi_t|^2dQ\\
&+\int_s^t\int_M|\frac{\tilde{\nabla}\psi_t}{\partial t}|^2|\tilde{\nabla}\psi_t|^2dQ
+\int_s^t\int_M\big(|\frac{\partial\phi_t}{\partial t}|^2+|\frac{\tilde{\nabla}\psi_t}{\partial t}|^2\big)dQ\big).
\end{align*}
The last term can be bounded in terms of the initial data and the \(L^2\)-norm of \(\psi_t\) by
Lemma \ref{lemma-surface-estimates-initial-data}.
We use another type of Sobolev inequality (similar to (\ref{local-sobolev-inequality}) for \(|t-s|\leq 1\)) 
to bound the mixed terms like
$\int_s^t\int_M|\frac{\partial\phi_t}{\partial t}|^2|\tilde{\nabla}\psi_t|^2dQ$, more precisely
\begin{align*}
\int_s^t\int_M|d\phi_t|^2|\frac{\partial\phi_t}{\partial t}|^2dQ\leq&
\big(\int_s^t\int_M|d\phi_t|^4dQ\big)^\frac{1}{2}
\big(\sup_{s\leq\theta\leq t} \int_M|\frac{\partial\phi}{\partial t}(\cdot,\theta)|^2 dM 
+\int_s^t\int_M |\nabla\frac{\partial\phi_t}{\partial t}|^2dQ\big)
\end{align*}
and similarly for both of the other two terms.

Choosing $t-s<\delta_2$ sufficiently small, applying the Sobolev inequality and
the estimates from Corollary \ref{surface-corollay-small-l4}, we can absorb part of the
right hand side in the left and obtain
\begin{align*}
\int_M\big(|\frac{\partial\phi_t(\cdot,t)}{\partial t}|^2&+|\frac{\tilde{\nabla}\psi_t(\cdot,t)}{\partial t}|^2\big)dM 
\leq \inf_{t-\delta_2\leq s\leq t}C\int_M\big(|\frac{\partial\phi_t(\cdot,s)}{\partial t}|^2+|\frac{\tilde{\nabla}\psi_t(\cdot,s)}{\partial t}|^2\big)dM+C.
\end{align*}
Finally, we estimate the infimum by the mean value, more precisely
\begin{align*}
\sup_{2\tau\leq t\leq T_1}\int_M\big(|\frac{\partial\phi(\cdot,t)}{\partial t}|^2+|\frac{\tilde{\nabla}\psi(\cdot,t)}{\partial t}|^2\big)dM
&\leq C(1+\tau^{-1})\int_s^t\int_M\big(|\frac{\partial\phi_t}{\partial t}|^2+|\frac{\tilde{\nabla}\psi_t}{\partial t}|^2\big)dQ+C\\
&\leq C(1+\tau^{-1}).
\end{align*}
Hence, we get the desired bound.
\end{proof}

\begin{Cor}
\label{theorem-surface-bounding-second-derivatives}
Let \((\phi_t,\psi_t)\in V\) be a solution of (\ref{evolution-phi-surface}) and (\ref{evolution-psi-surface}).
Assume \(|\psi_t|_{L^\infty(M\times [0,T_1))}\leq C\) and 
\(\sup_{(x,t)\in M\times[0,T_1)}F(\phi_t,\psi_t,B_R(x))<\delta_1\). Then we have
\begin{equation}
\int_M(|\nabla^2\phi(\cdot,t)|^2+\epsilon^2|\tilde{\nabla}^2\psi(\cdot,t)|^2)dM\leq C,
\end{equation}
where the constant \(C\) depends on $M,N,R,\delta_1,\delta_2,\epsilon,\tau,\psi_0,d\phi_0\) and \(\tilde{\nabla}\psi_0$.
\end{Cor}

\begin{proof}
With the help of the previous estimates we can now bound the full second derivatives of \((\phi_t,\psi_t)\) in \(L^2\).
By the evolution equations (\ref{evolution-phi-surface}), (\ref{evolution-psi-surface}) and Young's inequality,
we find
\begin{align*}
\int_M|\tau(\phi_t)|^2dM\leq&C\int_M(|\cR(\phi_t,\psi_t)|^2+\epsilon^2|\cR_c(\phi_t,\psi_t)|^2+|\frac{\partial\phi_t}{\partial t}|^2)dM\\
\leq& C\int_M(|\psi_t|^4|d\phi_t|^2+\epsilon^2|\psi_t|^2|\tilde{\nabla}\psi_t|^2|d\phi_t|^2+|\frac{\partial\phi_t}{\partial t}|^2)dM,\\
\int_M|\tilde{\nabla}_{e_\alpha}^*\tilde{\nabla}_{e_\alpha}\psi_t|^2dM\leq& C\int_M(|\D\psi_t|^2+|\frac{\tilde{\nabla}\psi_t}{\partial t}|^2)dM.
\end{align*}

The assertion then follows from applying the local Sobolev inequality (\ref{local-sobolev-inequality}) and the 
Bochner formulas \eqref{lemma-surface-bochner-phi}, \eqref{lemma-surface-bochner-psi} with \(\delta_1\) small enough.
\end{proof}

\begin{Cor}[Higher regularity]
Suppose the pair \((\phi_t,\psi_t)\) is a weak solution of (\ref{evolution-phi-surface}) and (\ref{evolution-psi-surface}).
The pair \((\phi_t,\psi_t)\) is smooth as long as \(\delta_1,\delta_2\) are small enough.
\end{Cor}
\begin{proof}
Since we have a bound on the \(L^2\)-norm of the second derivatives 
of \(\phi_t\) and \(\psi_t\) by (\ref{theorem-surface-bounding-second-derivatives}),
we can apply the Sobolev embedding theorem and get that both
\(|d\phi_t|\in L^p\) and \(|\tilde{\nabla}\psi_t|\in L^p\) for \(p<\infty\). From the evolution equations
(\ref{evolution-phi-surface}) and (\ref{evolution-psi-surface})
we may conclude that \(\big|\frac{\partial\phi_t}{\partial t}\big|,|\nabla^2\phi_t|\in L^p\)
and also \(\big|\frac{\tilde{\nabla}\psi_t}{\partial t}\big|,|\tilde{\nabla}^2\psi_t|\in L^p\).
By the regularity theory for parabolic partial differential equations we obtain that
\(|d\phi_t|\) and \(|\tilde{\nabla}\psi_t|\) are Hölder continuous,
see \cite{MR0241822}, Theorem IV.9.1 and Lemma II.3.3.
At this point the smoothness of the pair \((\phi,\psi)\) follows from a standard bootstrap argument using
Schauder theory, for more details see Theorem 3.24 in \cite{phd}.
\end{proof}

\section{Long-time Existence and Singularities}
In this section we establish the existence of a long-time solution to the evolution equations.
Thus, we first of all derive a uniqueness and stability result. To avoid the problem of identifying sections in different vector bundles, 
we will make use of the Nash embedding theorem. 
Hence, assume \(N\subset\R^q\) isometrically and denote the isometric embedding by \(\iota\).
Then, \(u=\iota\circ\phi\colon M\to \R^q\) can be thought of as a vector-valued function.
The vector spinor \(\psi\) turns into a vector of usual spinors \(\psi=(\psi^1,\ldots,\psi^q)\)
with \(\psi^i\in\Gamma(\Sigma M),~i=1,\ldots,q\). The condition that \(\psi\) is along the map \(\phi\) is encoded by
\[
\sum_{i=1}^q\nu_i\psi^i=0 \qquad \textrm{for a normal vector } \nu\in\R^q \textrm{ at } \phi(x).
\]
Now, the function \(u\) satisfies the following equation:
\begin{align}
\label{evolution-u-Rq}
\nonumber
\big(\frac{\partial}{\partial t}-\Delta\big)u=&-\sff_u(du,du)
-P(\sff(du(e_\alpha),e_\alpha\cdot\psi),\psi)-\epsilon B(du,\psi,du,\psi) \\
&+\epsilon P(\sff(du(e_\alpha),\psi),\nabla_{e_\alpha}\psi)
 -\epsilon P(\sff(du(e_\alpha),\nabla_{e_\alpha}\psi),\psi)
\end{align}
with the initial condition \(u_0=\iota(\phi_0)\) and 
\begin{align*}
B_u(du,\psi,du,\psi)=\frac{\partial u^i}{\partial x_\alpha}\frac{\partial u^k}{\partial x_\alpha}\langle\psi^{l},\psi^{j}\rangle 
\big(P(\sff_u(\partial_{y^i},\partial_{y^m}),\partial_{y^j})\Gamma^m_{kl}
-P(\sff_u(\partial_{y^i},\partial_{y^j}),\partial_{y^m})\Gamma^m_{kl}\big).
\end{align*}
For the spinor \(\psi\in\Gamma(\Sigma M\otimes T\R^q)\), we get the following evolution equation
\begin{align}
\label{evolution-psi-Rq}
\nonumber
\big(\frac{\nabla}{\partial t}-\epsilon\Delta\big)\psi=&-\p\psi
+\sff(du(e_\alpha),e_\alpha\cdot\psi)+\sff(\frac{\partial u}{\partial t},\psi)
-2\epsilon\sff(du(e_\alpha),\tilde{\nabla}_{e_\alpha}\psi)\\
 &-\epsilon(\nabla_{e_\alpha}\sff)(du(e_\alpha),\psi))-\epsilon\sff(\tau(u),\psi)
\end{align}
with the initial condition \(\psi_0=d\iota(\psi'_0)\), where \(\psi'_0\in\Gamma(\Sigma M\otimes\phi_0^{-1}TN)\).
For a derivation of \eqref{evolution-u-Rq} and \eqref{evolution-psi-Rq} see \cite{phd}, Section 3.4. 
Here, \(\sff\) is the second fundamental form of the embedding and
\(P\) denotes the shape operator. By projecting to a tubular neighborhood \(\tilde{N}\) of \(\iota(N)\subset\R^q\)
we can think of \(\sff\) as a vector-valued function in \(\R^q\). For more details, see \cite{MR1896863}, p. 132.

Assuming \(|\psi_t|_{L^\infty(M\times [0,T))}\leq C\) and \(\sup_{(x,t)\in M\times[0,T)}F(\phi_t,\psi_t,B_R(x))<\delta_1\)
we obtain by Corollary \ref{theorem-surface-bounding-second-derivatives} and the Sobolev embedding theorem that
\begin{equation}
\label{l4-bound-du}
\int_M|du|^4dM \leq C,\qquad \int_M|\nabla\psi|^4dM\leq C
\end{equation}
such that we can prove the following
\begin{Prop}[Stability and Uniqueness]
\label{surface-prop-uniqueness}
Let \((\phi,\psi)\in V\) and \((\xi,\chi)\in V\) be
solutions of (\ref{evolution-phi-surface}) and (\ref{evolution-psi-surface}),
where \(\psi\in\Gamma(\Sigma M\otimes\phi^{-1}TN)\) and moreover \(\chi\in\Gamma(\Sigma M\otimes\xi^{-1}TN)\).
In addition, suppose that \(|\psi|_{L^\infty(M\times [0,T))}\leq C\) and \(|\chi|_{L^\infty(M\times [0,T))}\leq C\).
If the initial data coincide, \((\phi_0,\psi_0)=(\xi_0,\chi_0)\), 
then we have $(\phi,\psi)=(\xi,\chi)$ throughout $M\times[0,T)$.
\end{Prop}
\begin{proof}
We follow \cite{MR2431434}, p.\ 235.
We regard \(u,v\) as vector-valued functions in \(\R^q\) with \(u=\iota\circ\phi,v=\iota\circ\xi\). 
The spinors \(\psi\) and \(\chi\) are defined along the maps \(u\) and \(v\).
We set 
\[
h(x,t)=(h_1(x,t),h_2(x,t))=(u(x,t)-v(x,t),\psi(x,t)-\chi(x,t)).
\]
First, we study the evolution of \(h_1\) and \(h_2\) separately and add up both contributions in the end.
We compute using (\ref{evolution-u-Rq})
\begin{align*}
\frac{\partial}{\partial t}\frac{1}{2}\int_M &|h_1|^2dM=-\int_M|dh_1|^2dM
- \int_M\langle\sff_u(du,du)-\sff_v(dv,dv),h_1\rangle dM\\
&-\int_M\langle h_1,P(\sff_u(du(e_\alpha),e_\alpha\cdot\psi),\psi)-P(\sff_v(dv(e_\alpha),e_\alpha\cdot\chi),\chi)\rangle dM\\
&-\epsilon\int_M\langle h_1,P(\sff_u(du(e_\alpha),\nabla_{e_\alpha}\psi),\psi)-P(\sff_v(dv(e_\alpha),\nabla_{e_\alpha}\chi),\chi)\rangle dM\\
&+\epsilon\int_M\langle h_1,P(\sff_u(du(e_\alpha),\psi),\nabla_{e_\alpha}\psi)-P(\sff_v(dv(e_\alpha),\chi),\nabla_{e_\alpha}\chi)\rangle dM\\
&-\epsilon\int_M\langle h_1,B_u(du,\psi,du,\psi)-B_v(dv,\chi,dv,\chi)\rangle dM.
\end{align*}
We estimate the right hand side in terms of \(h_1\) and \(h_2\), where we apply the pointwise bounds on \(\psi\) and \(\chi\).
We will sketch this in detail for the term with the second fundamental forms, the other terms can then be treated similarly.
Rearranging
\[
\sff_u(du,du)-\sff_v(dv,dv)=
(\sff_u-\sff_v)(du,du)+\sff_v(du-dv,du)+\sff_v(dv,du-dv)
\]
and applying the mean value theorem, we find
\[
|\langle \sff_u(du,du)-\sff_v(dv,dv),u-v\rangle|\leq C(|du|^2|u-v|^2+(|du|+|dv|)|h_1||dh_1|).
\]
Hence, we obtain the bound
\begin{align*}
|\langle \sff_u(du,du)-\sff_v(dv,dv),h_1\rangle|
&\leq C(|h_1|^2(|du|^2+|dv|^2)+\frac{1}{8}|dh_1|^2.
\end{align*}
Using the pointwise bound on the spinors we find
\begin{align*}
|\langle h_1,P(\sff_u(du(e_\alpha),e_\alpha\cdot\psi)&,\psi)-P(\sff_v(dv(e_\alpha),e_\alpha\cdot\chi),\chi)\rangle| \\
&\leq C(|du|^2|h_1|^2+|h_1|^2+|dv|^2|h_2|^2)+\frac{1}{8}|dh_1|^2
\end{align*}
and
\begin{align}
\label{uniqueness-estimate-1}
|\langle h_1,&P(\sff_u(du(e_\alpha),\psi),\nabla_{e_\alpha}\psi)-P(\sff_v(dv(e_\alpha),\chi),\nabla_{e_\alpha}\chi)\rangle| \\
\nonumber &\leq C(|h_1|^2(|du|^2+|\nabla\psi|^2+|dv|^2)+|h_2|^2|\nabla\psi|^2)+\frac{1}{8}|dh_1|^2+\frac{1}{8}|\nabla h_2|^2.
\end{align}
Note that the contribution 
\[
\epsilon\langle h_1,P(\sff_u(du(e_\alpha),\nabla_{e_\alpha}\psi),\psi)-P(\sff_v(dv(e_\alpha),\nabla_{e_\alpha}\chi),\chi)\rangle\\
\]
can be estimated the same way as \eqref{uniqueness-estimate-1}. In addition, we have
\begin{align*}
|\langle h_1,B_u(du,\psi,du,\psi)&-B_v(dv,\chi,dv,\chi)\rangle|
\leq C(|h_1|^2(|du|^2+|dv|^2)+|h_2|^2|dv|^2)+\frac{1}{8}|dh_1|^2.
\end{align*}

We now turn to the function \(h_2\). With the help of (\ref{evolution-psi-Rq}) we find 
\begin{align}
\label{evolution-norm-h_2}
\frac{\partial}{\partial t}\frac{1}{2}\int_M|h_2|^2dM=&-\int_M\langle\p h_2,h_2\rangle dM-\epsilon\int_M|\nabla h_2|^2dM\\
\nonumber &-\epsilon\int_M\langle h_2,(\nabla_{e_\alpha}\sff_u)(du(e_\alpha),\psi)-(\nabla_{e_\alpha}\sff_v)(dv(e_\alpha),\chi)\rangle dM.
\end{align}
The other terms involving the second fundamental form vanish since \(\sff\perp\psi\) and
we may estimate
\[
\langle\p h_2,h_2\rangle\leq C|h_2|^2+\frac{\epsilon}{8}|\nabla h_2|^2.
\]

To estimate the last term in \eqref{evolution-norm-h_2}, we rearrange and estimate
\begin{align*}
|\langle \psi-\chi,(\nabla_{e_\alpha}\sff_u)(du(e_\alpha),\psi)&-(\nabla_{e_\alpha}\sff_v)(dv(e_\alpha),\chi)\rangle|
\leq C(|h_1|^2|du|^2+|dv|^2|h_2|^2)+\frac{1}{8}|dh_1|^2.
\end{align*}

Adding up the inequalities for \(|h_1|^2\) and \(|h_2|^2\) and applying the 
Sobolev embedding theorem, we find
\begin{align*}
&\frac{\partial}{\partial t}\frac{1}{2}\int_M(|h_1|^2+|h_2|^2)dM+\frac{1}{2}\int_M(|dh_1|^2+\epsilon|\nabla h_2|^2)dM \\
&\leq C\int_M(|h_1|^2+|h_2|^2)dM+C\int_M(|h_1|^2+|h_2|^2)(|du|^2+|dv|^2+|\nabla\psi|^2)dM.
\end{align*}
Applying \eqref{l4-bound-du} and using the Sobolev embedding theorem
the last term on the right hand side can be estimated as
\begin{align*}
C\int_M|h_1|^2|du|^2dM\leq &C\big(\int_M|h_1|^4dM\big)^\frac{1}{2}\big(\int_M|du|^4dM\big)^\frac{1}{2} \\
\leq& C\big(\int_M|h_1|^2dM\big)^\frac{1}{2}\big(\int_M|dh_1|^2dM\big)^\frac{1}{2} \\
\leq&\frac{1}{8}\int_M|dh_1|^2dM+C\int_M|h_1|^2dM
\end{align*}
and the other contributions can be treated similarly. Hence, we find
\begin{align*}
\frac{\partial}{\partial t}\frac{1}{2}&\int_M(|h_1|^2+|h_2|^2)dM\leq C\int_M(|h_1|^2+|h_2|^2)dM.
\end{align*}
Integrating with respect to \(t\) and using that \(h_1(0)=h_2(0)=0\) we may follow that
\(h_1=h_2=0\) for all \(t\in [0,T)\), which proves the claim.
\end{proof}

\begin{Prop}[Long-time Existence]
Let \((\phi_t,\psi_t)\in V\) be a solution of (\ref{evolution-phi-surface}) and (\ref{evolution-psi-surface}).
Assume that \(|\psi_t|_{L^\infty(M\times [0,T))}\leq C\).
Then the evolution equations admit a unique weak solution for $0\leq t<\infty$.
\end{Prop}
\begin{proof}
The first singular time \(T_0\) is characterized by the condition
\[
\limsup_{t\to T_0}F(\phi_t,\psi_t,B_R(x))\geq \delta_1.
\]
Since we have $\partial_t\phi,\tilde{\nabla}_t\psi\in L^2(M\times [0,T_0))$ and
also $F(\phi_t,\psi_t)\leq CF(\phi_0,\psi_0)+C$ for \mbox{$0<t<T_0$},
there exists 
\[(\phi(\cdot,T_0),\psi(\cdot,T_0))\in H^{1}(M,N)\times H^{1}(M,\Sigma M\otimes\phi_t^{-1}TN)
\]
such that 
\[(\phi(\cdot,t),\psi(\cdot,t))\to (\phi(\cdot,T_0),\psi(\cdot,T_0))\]
weakly in \(H^{1}(M,N)\times H^{1}(M,\Sigma M\otimes\phi_t^{-1}TN)$ as $t$ approaches $T_0$. 
In particular, we have
\[
F(\phi_{T_0},\psi_{T_0})\leq\liminf_{s\to t} CF(\phi_s,\psi_s) +C\leq CF(\phi_t,\psi_t)+C,\qquad 0\leq t\leq T_0.
\]
Now let
\((\tilde{\phi}_t,\tilde{\psi}_t)\colon (M\times[T_0,T_0+T_1)\to N)\times (M\times[T_0,T_0+T_1)\to\Sigma M\otimes\phi_t^{-1}TN)\)
be a solution of (\ref{evolution-phi-surface}) and (\ref{evolution-psi-surface}).
Assume that $(\tilde{\phi},\tilde{\psi})(x,t)=(\phi,\psi)(x,t)$.
We define
\[
(\hat{\phi}_t,\hat{\psi}_t)=
\begin{cases}
(\phi_t,\psi_t), &\qquad 0\leq t\leq T_0,\\
(\tilde{\phi}_t,\tilde{\psi}_t), &\qquad T_0\leq t\leq T_0+T_1.
\end{cases}
\]
Now
\(
(\hat{\phi}_t,\hat{\psi}_t)\colon (M\times[0,T_0+T_1)\to N)\times (M\times[0,T_0+T_1)\to\Sigma M\times\hat{\phi}_t^{-1}TN)
\)
is a weak solution of (\ref{evolution-phi-surface}) and (\ref{evolution-psi-surface}).
By iteration, we obtain a weak solution $(\phi_t,\psi_t)$ on a maximal time interval $T_0+\delta$
for some $\delta>0$. If $T_0+\delta<\infty$ then by the above argument the solution $(\phi_t,\psi_t)$
may be extended to infinity, hence $T_0+\delta=\infty$. The uniqueness follows from Proposition (\ref{surface-prop-uniqueness}).
\end{proof}

\begin{Prop}
Assume \((\phi_t,\psi_t)\) is
a solution of (\ref{evolution-phi-surface}) and (\ref{evolution-psi-surface})
satisfying \eqref{condition-poincare}.
There are only finitely many singular points \((x_k,t_k),1\leq k\leq K\).
The number \(K\) depends on \(M,\epsilon,\psi_0,d\phi_0\) and \(\tilde{\nabla}\psi_0\).
\end{Prop}
\begin{proof}
We follow the presentation in \cite{MR2431658}, p.\ 138, for the harmonic map heat flow.
We assume that $T_0>0$ is the first singular time and define the singular set as
\begin{equation}
\label{blowup-set}
S(\phi,\psi,T_0)=\bigcap_{R>0}\big\{
x\in M\mid\limsup_{t\to T_0}F(\phi_t,\psi_t,B_R(x))\geq \delta_1
\big\}.
\end{equation}
Now, let $\{x_j\}^K_{j=1}$ be any finite subset of $S(\phi,\psi,T_0)$.
Then we have for $R>0$
\[
\limsup_{t\to T_0}\int_{B_R(x_j)}(|d\phi|^2+\epsilon|\tilde{\nabla}\psi|^2)dM\geq\delta_1,\qquad 1\leq j\leq K.
\]
By (\ref{surface-inequality-F-local}) we have the following local inequality for the quantity \(F(\phi_t,\psi_t,B_R)\)
\begin{equation}
\label{monotonicity-F}
F(\phi_t,\psi_t,B_R(x))\leq 2E_\epsilon(\phi_0,\psi_0,B_{2R}(x))+\delta_3\frac{T}{R^2}+\frac{1}{\epsilon}\int_{B_R}|\psi_t|^2dM
\end{equation}
with \(\delta_3=C\int_M(|d\phi_t|^2+\epsilon^2|\tilde{\nabla}\psi_t|^2+|\psi_t|^2)dM\).
Since \(E_\epsilon(\phi_t,\psi_t)\leq E_\epsilon(\phi_0,\psi_0)\) we obtain
\begin{align*}
-\frac{1}{\sqrt{2}}\big(\int_M|\psi_t|^2dM\big)^\frac{1}{2}\big(\int_M|\tilde{\nabla}\psi_t|^2dM\big)^\frac{1}{2}
&+F(\phi_t,\psi_t) \\
&\leq F(\phi_0,\psi_0)+\frac{1}{\sqrt{2}}\big(\int_M|\psi_0|^2dM\big)^\frac{1}{2}\big(\int_M|\tilde{\nabla}\psi_0|^2dM\big)^\frac{1}{2}.
\end{align*}
Recall that by assumption \eqref{condition-poincare} we have, 
\[
\int_M|\psi_t|^2dM\leq\delta_5\int|\tilde{\nabla}\psi|^2dM,
\]
where we have renamed the positive constant \(c_1\) to \(\delta_5\).
From this we obtain the global estimate (with \(\epsilon\) suitably large)
\[
F(\phi_t,\psi_t)\leq\delta_4 F(\phi_0,\psi_0)
\]
for a positive constant \(\delta_4=\frac{\max\{\frac{\sqrt{2}}{\delta_5}+\epsilon,1\}}{\min\{\frac{-\sqrt{2}}{\delta_5}+\epsilon,1\}}\).
We choose $R>0$ such that all the $B_{2R}(x_j),1\leq j\leq K$ are mutually disjoint
and small enough to have
\[
\frac{1}{\epsilon}\int_{B_R}|\psi_t|^2dM\leq\frac{\delta_1}{4}.
\]
Then, we have by (\ref{monotonicity-F})
\begin{align*}
K\delta_1\leq& \sum_{j=1}^K\limsup_{t\to T_0} F(\phi_t,\psi_t,B_R(x_j)) \\
\leq& \sum_{j=1}^K\big(\limsup_{t\to T_0} 2E_\epsilon(\phi_\tau,\psi_\tau,B_{2R}(x_j))+\frac{\delta_1}{2}\big) \\
\leq& 2E_\epsilon(\phi_\tau,\psi_\tau)+\frac{K\delta_1}{2} \\
\leq& 2E_\epsilon(\phi_0,\psi_0)+\frac{K\delta_1}{2}
\end{align*}
for any $\tau\in[T_0-\frac{\delta_1R^2}{4\delta_3},T_0]$.
We conclude that
\[
K\leq 4\frac{E_\epsilon(\phi_0,\psi_0)}{\delta_1},
\]
which implies the finiteness of the singular set $S(\phi,\psi,T_0)$.
Our next aim is to show that there are only finitely many singular spatial points.
Therefore we set
\[
\tilde{M}=M\setminus\bigcup_{1\leq j\leq K}B_{2R}(x_j)
\]
and in addition, we calculate
\begin{align}
\label{inequality-F-singularities}
F(\phi_{T_0},\psi_{T_0})=&\lim_{R\to 0} F(\phi_{T_0},\psi_{T_0},\tilde{M})\\
\nonumber\leq&\lim_{R\to 0}\limsup_{t\to T_0} F(\phi_{t},\psi_{t},\tilde{M})\\
\nonumber =&F(\phi_t,\psi_t)-\lim_{R\to 0}\sum_{j=1}^K\liminf_{t\to T_0} F(\phi_t,\psi_t,B_{2R}(x_j))\\
\nonumber\leq&\delta_4 F(\phi_0,\psi_0)-\lim_{R\to 0}\sum_{j=1}^K\limsup_{t\to T_0} F(\phi_t,\psi_t,B_{R}(x_j))\\
\nonumber\leq&\delta_4 F(\phi_0,\psi_0)-K\delta_1.
\end{align}
Now suppose $T_0<\ldots<T_j$ are $j$ singular times and by $K_0,\ldots,K_j$
we denote the number of singular points at each singular time.
Set
\[
(\phi_i,\psi_i)=\lim_{t\to T_i}(\phi_t,\psi_t),\qquad 0\leq i\leq j.
\]
By iterating (\ref{inequality-F-singularities}) we get
\begin{align*}
F(\phi_j,\psi_j)\leq&\delta_4 F(\phi_{j-1},\psi_{j-1})-\delta_1 K_{j-1} \\
\leq& \delta_4^2 F(\phi_{j-2},\psi_{j-2})-\delta_1( K_{j-1}+\delta_4 K_{j-2})\\
\leq&\ldots \\
\leq&\delta_4^jF(\phi_0,\psi_0)-\delta_1\sum_{i=0}^{j-1}K_i\delta_4^{j-i-1},
\end{align*}
which can be rearranged as
\begin{equation}
\label{finitely-many-spatial-singularities}
\sum_{i=0}^{j-1}K_i\delta_4^{-i-1}\leq \frac{F(\phi_0,\psi_0)}{\delta_1}.
\end{equation}
We conclude that there are only finitely many singularities.
\end{proof}
\begin{Bem}
If we compare the bound on the number of singularities of the regularized Dirac-harmonic map heat flow
with the bound on the number of singularities in the harmonic map heat flow, then we realize
that the former can encounter more singularities.
In the case of the harmonic map heat flow we would have \(\delta_4=1\) and \(F(\phi_0,\psi_0)=\frac{1}{2}\int_M|d\phi_0|^2\),
which lowers the upper bound in (\ref{finitely-many-spatial-singularities}).
\end{Bem}

\section{Convergence and Blowup Analysis}
In this section we discuss the convergence of the evolution equations (\ref{evolution-phi-surface}) and (\ref{evolution-psi-surface}).
In addition, we address the problem of blowing up the singular points.
% \begin{Prop}
% \label{theorem-surface-convergence}
% Let \((\phi_t,\psi_t)\in V\) be a solution of (\ref{evolution-phi-surface}) and (\ref{evolution-psi-surface}).
% Moreover, assume that \(|\psi_t|_{L^\infty(M\times [0,\infty))}\leq C\).
% Then the pair \((\phi_t,\psi_t)\) converges strongly in \(L^2\) to a regularized 
% Dirac-harmonic map on \(M\setminus\{x_1,\ldots,x_k\}\).
% The limiting map \((\phi_\infty,\psi_\infty)\) is smooth on \(M\setminus\{x_1,\ldots,x_k\}\).
% \end{Prop}

\begin{Prop}
\label{theorem-surface-convergence}
Let \((\phi_t,\psi_t)\in V\) be a solution of (\ref{evolution-phi-surface}) and (\ref{evolution-psi-surface}).
Moreover, assume that \(|\psi_t|_{L^\infty(M\times [0,\infty))}\leq C\).
Then the pair \((\phi_t,\psi_t)\) converges weakly in $H^{1}(M,N)\times H^{1}(M,\Sigma M\otimes\phi_t^{-1}TN)$ 
and strongly in the space $W^{2,2}_{loc}(M\setminus\{x_k,t_k=\infty\},N)\times W^{2,2}_{loc}(M\setminus\{x_k,t_k=\infty\},\Sigma M\otimes\phi_t^{-1}TN)$ 
to a regularized Dirac-harmonic map. The limiting map \((\phi_\infty,\psi_\infty)\) is smooth on \(M\setminus\{x_1,\ldots,x_k\}\).
\end{Prop}

\begin{proof}
Since we have a uniform bound on the \(L^2\)-norm of the \(t\) derivatives of \((\phi_t,\psi_t)\) 
by Lemma \ref{lemma-surface-estimates-initial-data}, we can achieve for $t_m\to\infty$ suitably 
\[
\int_M\big(\big|\frac{\partial\phi_t}{\partial t}\big|^2+\big|\frac{\tilde{\nabla}\psi_t}{\partial t}\big|^2\big)dM\big|_{t=t_m}\to 0
\]
and in addition, we suppose that $T=\infty$ is non-singular
\[
\limsup_{t\to\infty} (\sup_{x\in M}F(\phi_t,\psi_t,B_R(x)))<\delta_1
\]
for some \(R>0\).
By (\ref{theorem-surface-bounding-second-derivatives}) we have a bound on the second derivatives
\[
\int_M\big(|\nabla^2\phi|^2(\cdot,t_m)+\epsilon^2|\tilde{\nabla}^2\psi|^2(\cdot,t_m)\big)dM\leq C
\]
and due to the Rellich-Kondrachov embedding theorem we may assume that 
\begin{align*}
\phi(\cdot,t_m)\to \phi_\infty& \qquad \textrm{strongly in}~W^{1,p}(M,N), \\
\psi(\cdot,t_m)\to \psi_\infty& \qquad \textrm{strongly in}~W^{1,p}(M,\Sigma M\otimes\phi_{t_m}^{-1}TN)
\end{align*}
for any $p<\infty$.
But then by (\ref{evolution-phi-surface}) and (\ref{evolution-psi-surface}) we get
convergence of the evolution equations
\begin{align}
\label{surface-solution-phi}\tau(\phi_\infty)=&\cR(\phi_\infty,\psi_\infty)+\epsilon\cR_c(\phi_\infty,\psi_\infty),\\
\label{surface-solution-psi}\epsilon\tilde{\Delta}\psi_\infty=&\D\psi_\infty
\end{align}
in $L^2$, the pair $(\phi_\infty,\psi_\infty)$ is a regularized Dirac-harmonic map,
which satisfies
\(
(\phi_\infty,\psi_\infty)\in W^{2,2}(M,N)\times W^{2,2}(M,\Sigma M\otimes\phi_\infty^{-1}TN).
\)
\\
If $T=\infty$ is singular, meaning that at the points
$\{x_1,\ldots,x_k\}$
\[
\limsup_{t\to\infty} F(\phi_t,\psi_t,B_R(x_j))\geq\delta_1,\qquad 1\leq j\leq k
\]
for all $R>0$, then for suitable numbers $t_m\to\infty$ the family $(\phi_{t_m},\psi_{t_m})$
will be bounded in $W^{2,2}_{loc}(M,N)\times W^{2,2}_{loc}(M,\Sigma M\otimes\phi_{t_m}^{-1}TN)$ on the set \(M\setminus\{x_1,\ldots,x_k\}\).
Consequently, the family \((\phi_{t_m},\psi_{t_m})\) will accumulate as follows
\begin{eqnarray*}
&&\phi_\infty\colon M\setminus\{x_1,\ldots,x_k\}\to N, \\
&&\psi_\infty\colon M\setminus\{x_1,\ldots,x_k\}\to \Sigma (M\setminus\{x_1,\ldots,x_k\})\otimes\phi_{\infty}^{-1} TN.
\end{eqnarray*}
We set \(\tilde{M}:=M\setminus\{x_1,\ldots,x_k\}\).
Concerning the regularity of \((\phi_\infty,\psi_\infty)\) on \(\tilde{M}\), 
we have \(\phi_\infty\in W^{1,p}_{loc}(\tilde{M},N)\) for any \(0<p<\infty\), since \(\phi_\infty\in W^{2,2}_{loc}(\tilde{M},N)\).
In addition, we have \(\psi_\infty\in W^{2,2}_{loc}(\tilde{M},\Sigma\tilde{M}\otimes\phi_\infty^{-1}TN)\) and consequently also
\(\psi_\infty\in W^{1,p}_{loc}(\tilde{M},\Sigma\tilde{M}\otimes\phi_\infty^{-1}TN)\) for any \(0<p<\infty\).
Hence, the right hand sides of both (\ref{surface-solution-phi}) and (\ref{surface-solution-psi}) are in \(L^p_{loc}\) for \(2<p<\infty\).
Writing \(\tau(\phi)=\Delta\phi+\Gamma(\phi)(d\phi,d\phi)\) and by elliptic estimates for second order operators 
we then get \(\phi_\infty\in W^{2,p}_{loc}(\tilde{M},N)\) for any \(0<p<\infty\).
The smoothness of \((\phi_\infty,\psi_\infty)\) then follows from a standard bootstrap argument.
\end{proof}
This completes the proof of Theorem \ref{theorem-surface}.
\par\medskip
Our next aim is to get a better understanding of the singular points \((x_k,t_k)\).
In the case of the harmonic map heat flow one can perform a blowup analysis,
which finally leads to the ``bubbling off of harmonic spheres'', see for example \cite{MR826871}.
The important ingredient in that calculation is the fact that one can perform a parabolic rescaling
of the evolution equation for harmonic maps. 
Thus, let us analyze the scaling of the regularized Dirac-harmonic heat flow.
\begin{Bem}
By regularizing the functional $E(\phi,\psi)$, we haven broken
the conformal invariance and consequently the evolution equations for $(\phi_t,\psi_t)$
do not scale in a ``nice'' way.
Nevertheless, it is possible to do a rescaling if one allows to rescale $\epsilon$ as well.
It is easy to see that the evolution equations \eqref{evolution-phi-surface} and \eqref{evolution-psi-surface}
are invariant under the following rescaling
\begin{align}
\phi(x,t)\to&\phi(x_0+Rx,t_0+R^2t), \\
\nonumber\psi(x,t)\to&\sqrt{R}\psi(x_0+Rx,t_0+Rt),\\
\nonumber\epsilon\to&\frac{\epsilon}{R}
\end{align}
for \(R>0\).
A dimensional analysis of the evolution equation for \(\psi\) also motivates to rescale \(\epsilon\).
Note that the two evolution equations scale differently. The evolution equation
for \(\phi\) scales like a heat type equation, whereas the evolution equation for \(\psi\)
scales like a first order evolution equation. However, it seems impossible to justify the rescaling of  \(\epsilon\) 
at a rigorous level.
\end{Bem}

\begin{Bem}
When analyzing the bubbling of Dirac-harmonic maps, it is important to have
control over the energy of the bubbles, such that now concentration phenomena can happen.
This control is usually given by what is called \emph{energy identity}.
For Dirac-harmonic maps the energy identity was established in \cite{MR2267756}, p.\ 131.

\begin{Dfn}
Let $(\phi_k,\psi_k):M\to N$ be a sequence of smooth Dirac-harmonic maps with uniformly bounded energy
\[
\int_M(|d\phi_k|^2+|\psi_k|^4)dM\leq C
\]
and furthermore assume that $(\phi_k,\psi_k)$ converges weakly to a Dirac-harmonic map $(\phi,\psi)$ in
$H^{1}(M,N)\times L^4(\Sigma M\otimes T\R^q)$. Then we call
\[
S:=\bigcap_{R>0}\{x\in M\mid\liminf_{k\to\infty}\int_{B_R(x)}(|d\phi_k|^2+|\psi_k|^4)dM>\delta\}
\]
the blow-up set of $\{\phi_k,\psi_k\}$.
\end{Dfn}
Note that the blow-up set for Dirac-harmonic maps differs from the blow-up set for 
regularized Dirac-harmonic maps \eqref{blowup-set} that we encountered when studying
the evolution equations.
\end{Bem}

\section{Removing the Regularization}
In this section we analyze the limit \(\epsilon\to 0\).
We have seen that the regularized Dirac-harmonic map heat flow
converges to a smooth regularized Dirac-harmonic map \((\phi_\infty,\psi_\infty)\) on \(M\) away from finitely many 
singular points. The smoothness of the limiting map depends on the estimates that were derived before.
Therefore the question is, which of these estimates we still need to control 
after taking the limit \(\epsilon\to 0\).
In particular, we would like to
\begin{enumerate}
\item Keep the number of singularities bounded,
\item Remove the singularities of the solution \((\phi_\infty,\psi_\infty)\),
\item Control the regularity of the solution \((\phi_\infty,\psi_\infty)\).
\end{enumerate}
Note that there is no preferred order in which these steps should be performed.
However, we cannot expect that the limit \(\epsilon\to 0\) will exist in general.
\begin{Bsp}
\begin{enumerate}
\item Assume that \(M=S^2\) and \(N=T^2\). In this case the Euler-Lagrange equations decouple
and we have to look for harmonic spinors on \(S^2\). It is well-known that these do
not exist \cite{MR1162671}. Consequently, the limit \(\epsilon\to 0\) cannot exist in this case and this fact should be reflected by the calculation.
\item If both \(M=N=T^2\), the Euler-Lagrange equations also decouple and we have
to look for harmonic spinors on \(T^2\). The two-dimensional torus has four spin structures
and not all of them admit harmonic spinors. Hence, the limit \(\epsilon\to 0\) cannot be trivial in this case, too.
\end{enumerate}
\end{Bsp}

\begin{Bem}
Both examples show that our approach using the \(L^2\)-gradient flow of the regularized functional \(E_\epsilon(\phi,\psi)\)
cannot detect the structures associated to the spinor bundle like the spin structure.
\end{Bem}

\subsection*{Number of singularities after \texorpdfstring{\(\epsilon\to 0\)}{}}
To study the dependence of the bound on the number of singularities on \(\epsilon\), 
let us analyze how the bound (\ref{finitely-many-spatial-singularities}) depends on \(\epsilon\).
Rearranging (\ref{finitely-many-spatial-singularities}) yields
\begin{equation}
\label{surface-number-singularities-epsilontozero}
\sum_{i=0}^{j-1}K_i\leq C\frac{F(\phi_0,\psi_0)(\epsilon)}{\delta_1(\epsilon)}.
\end{equation}
It is easy to see that 
\[
\lim_{\epsilon\to 0} F(\phi_0,\psi_0)=E(\phi_0)\leq C,
\]
but on the other hand the limit
\[
\lim_{\epsilon\to 0}\delta_1(\epsilon)
\]
does not exist in general as can easily be seen from the definition of \(\delta_1\).
Moreover, there is no cancellation of the different \(\epsilon\)'s on the right hand side of
(\ref{surface-number-singularities-epsilontozero}).

\subsection*{Removal of singularities after \texorpdfstring{\(\epsilon\to 0\)}{}}
To remove the singularities of the solution \((\phi_\infty,\psi_\infty)\) we would like 
to apply the following (Theorem 4.6 in \cite{MR2262709}, p.\ 426):
\begin{Satz}[Removable singularity theorem]
For $U\subset M$ let $(\phi,\psi)$ be a Dirac-harmonic map which is $C^\infty$ on $U\setminus\{p\}$ for some $p\in U$.
If 
\[
\int_U\big(|d\phi|^2+|\psi|^4\big)dM\leq C
\]
then $(\phi,\psi)$ extends to a $C^\infty$ solution on $U$.
\end{Satz}
In our case, the \(L^2\)-norm of \(d\phi_\infty\) can be bounded by plugging the spinor \(\psi_\infty\)
into the inequality for the energy functional \(E_\epsilon(\phi,\psi)\)
\[
\int_M|d\phi_\infty|^2dM\leq E_\epsilon(\phi_0,\psi_0).
\]
Unfortunately, we cannot  bound the \(L^4\)-norm of \(\psi_\infty\) after \(\epsilon\to 0\). 

\subsection*{Regularity of \texorpdfstring{\((\phi_\infty,\psi_\infty)\)}{} after \texorpdfstring{\(\epsilon\to 0\)}{}}
The regularity of Dirac-harmonic maps has been studied in \cite{MR2544729}.
\begin{Dfn}[Weakly Dirac-harmonic map]
A weak Dirac-harmonic map is a pair $(\phi,\psi)\in W^{1,2}(M,N)\times W^{1,\frac{4}{3}}{(M,\Sigma M\otimes\phi^{-1}TN})$,
which solves \eqref{euler-lagrange-phi-unregularized} and \eqref{euler-lagrange-psi-unregularized} in a weak sense.
\end{Dfn}
The relation between weak and smooth Dirac-harmonic maps in dimension two is given by the following (\cite{MR2544729}, Theorem 1.5, p.\ 3764)
\begin{Satz}
\label{regularity-surface}
Assume that \(M\) is a compact Riemannian spin surface and that the pair
$(\phi,\psi)\in W^{1,2}(M,N)\times W^{1,\frac{4}{3}}(M,\Sigma M\otimes\phi^{-1}TN)$ is a weak Dirac-harmonic map. 
Then the pair $(\phi,\psi)$ is smooth.
\end{Satz}
Hence, we have to ensure that the estimates necessary for the existence of a weakly Dirac-harmonic map
can be carried over to the limit \(\epsilon\to 0\).
Again, the regularity of the map \(\phi\) can be assured, but we do not have control
over \(\psi_\infty\) after  \(\epsilon\to 0\). 

\appendix
\section{}
The following Lemma combines the pointwise maximum principle 
with an integral norm. It can be thought of as a simple version of 
Moser's parabolic Harnack inequality.
\begin{Lem}
\label{maximum-principle-l2}
Assume that \((M,h)\) is a compact Riemannian manifold. If a function \(u(s,t)\geq 0\) satisfies
\[
\frac{\partial u}{\partial t}\leq \Delta u+Cu,
\]
and if in addition we have the bound
\[
U(t)=\int_Mu(s,t)dM\leq U_0,
\]
then there exists a uniform bound on 
\[
u(s,t)\leq e^CKU_0 
\]
with the constant \(K\) depending on \(M\).
\end{Lem}
\begin{proof}
A proof can for example be found in \cite{MR2744149}, p.\ 284.
\end{proof}

\emph{Acknowledgements:}
The author would like to thank the ``IMPRS for Geometric Analysis, Gravitation and String Theory'' 
for financial support.
In addition, the author gratefully acknowledges the support of the Austrian Science Fund (FWF) 
through the project P30749-N35 ``Geometric variational problems from string theory''.
\bibliographystyle{plain}
\bibliography{mybib}

\def\cprime{$'$}
\begin{thebibliography}{10}

\bibitem{MR3070562}
Bernd Ammann and Nicolas Ginoux.
\newblock Dirac-harmonic maps from index theory.
\newblock {\em Calc. Var. Partial Differential Equations}, 47(3-4):739--762,
  2013.

\bibitem{MR1162671}
Christian B{\"a}r.
\newblock Lower eigenvalue estimates for {D}irac operators.
\newblock {\em Math. Ann.}, 293(1):39--46, 1992.

\bibitem{MR1870959}
Michiel Bertsch, Roberta Dal~Passo, and Rein van~der Hout.
\newblock Nonuniqueness for the heat flow of harmonic maps on the disk.
\newblock {\em Arch. Ration. Mech. Anal.}, 161(2):93--112, 2002.

\bibitem{phd}
Volker Branding.
\newblock The evolution equations for {D}irac-harmonic maps, {P}h{D} thesis.
\newblock {\em http://opus.kobv.de/ubp/volltexte/2013/6420/}, 2013.

\bibitem{MR3305429}
Volker Branding.
\newblock Magnetic {D}irac-harmonic maps.
\newblock {\em Anal. Math. Phys.}, 5(1):23--37, 2015.

\bibitem{MR3333092}
Volker Branding.
\newblock Some aspects of {D}irac-harmonic maps with curvature term.
\newblock {\em Differential Geom. Appl.}, 40:1--13, 2015.

\bibitem{MR3493217}
Volker Branding.
\newblock Dirac-harmonic maps with torsion.
\newblock {\em Commun. Contemp. Math.}, 18(4):1550064, 19, 2016.

\bibitem{MR3558358}
Volker Branding.
\newblock Energy estimates for the supersymmetric nonlinear sigma model and
  applications.
\newblock {\em Potential Anal.}, 45(4):737--754, 2016.

\bibitem{MR3435758}
Volker Branding.
\newblock The evolution equations for regularized {D}irac-geodesics.
\newblock {\em J. Geom. Phys.}, 100:1--19, 2016.

\bibitem{MR3886921}
Volker Branding.
\newblock A vanishing result for the supersymmetric nonlinear sigma model in
  higher dimensions.
\newblock {\em J. Geom. Phys.}, 134:1--10, 2018.

\bibitem{MR3917346}
Volker Branding.
\newblock Energy methods for {D}irac-type equations in two-dimensional
  {M}inkowski space.
\newblock {\em Lett. Math. Phys.}, 109(2):295--325, 2019.

\bibitem{MR4034775}
Volker Branding.
\newblock Nonlinear {D}irac equations, monotonicity formulas and {L}iouville
  theorems.
\newblock {\em Comm. Math. Phys.}, 372(3):733--767, 2019.

\bibitem{MR3830277}
Volker Branding and Klaus Kr\"{o}ncke.
\newblock Global existence of {D}irac-wave maps with curvature term on
  expanding spacetimes.
\newblock {\em Calc. Var. Partial Differential Equations}, 57(5):Art. 119, 30,
  2018.

\bibitem{MR1180392}
Kung-Ching Chang, Wei~Yue Ding, and Rugang Ye.
\newblock Finite-time blow-up of the heat flow of harmonic maps from surfaces.
\newblock {\em J. Differential Geom.}, 36(2):507--515, 1992.

\bibitem{MR2370260}
Q.~Chen, J.~Jost, and G.~Wang.
\newblock Liouville theorems for {D}irac-harmonic maps.
\newblock {\em J. Math. Phys.}, 48(11):113517, 13, 2007.

\bibitem{MR2176464}
Qun Chen, J{\"u}rgen Jost, Jiayu Li, and Guofang Wang.
\newblock Regularity theorems and energy identities for {D}irac-harmonic maps.
\newblock {\em Math. Z.}, 251(1):61--84, 2005.

\bibitem{MR2262709}
Qun Chen, J{\"u}rgen Jost, Jiayu Li, and Guofang Wang.
\newblock Dirac-harmonic maps.
\newblock {\em Math. Z.}, 254(2):409--432, 2006.

\bibitem{MR3412386}
Qun Chen, J{\"u}rgen Jost, Linlin Sun, and Miaomiao Zhu.
\newblock Dirac-geodesics and their heat flows.
\newblock {\em Calc. Var. Partial Differential Equations}, 54(3):2615--2635,
  2015.

\bibitem{MR3044133}
Qun Chen, J{\"u}rgen Jost, and Guofang Wang.
\newblock The maximum principle and the {D}irichlet problem for
  {D}irac-harmonic maps.
\newblock {\em Calc. Var. Partial Differential Equations}, 47(1-2):87--116,
  2013.

\bibitem{MR3085099}
Qun Chen, J{\"u}rgen Jost, Guofang Wang, and Miaomiao Zhu.
\newblock The boundary value problem for {D}irac-harmonic maps.
\newblock {\em J. Eur. Math. Soc. (JEMS)}, 15(3):997--1031, 2013.

\bibitem{MR1701618}
Pierre Deligne, Pavel Etingof, Daniel~S. Freed, Lisa~C. Jeffrey, David Kazhdan,
  John~W. Morgan, David~R. Morrison, and Edward Witten, editors.
\newblock {\em Quantum fields and strings: a course for mathematicians. {V}ol.
  1, 2}.
\newblock American Mathematical Society, Providence, RI; Institute for Advanced
  Study (IAS), Princeton, NJ, 1999.
\newblock Material from the Special Year on Quantum Field Theory held at the
  Institute for Advanced Study, Princeton, NJ, 1996--1997.

\bibitem{MR0420708}
J.~Eells and J.~C. Wood.
\newblock Restrictions on harmonic maps of surfaces.
\newblock {\em Topology}, 15(3):263--266, 1976.

\bibitem{MR2138082}
Xiaoli Han.
\newblock Dirac-wave maps.
\newblock {\em Calc. Var. Partial Differential Equations}, 23(2):193--204,
  2005.

\bibitem{MR2860404}
Takeshi Isobe.
\newblock On the existence of nonlinear {D}irac-geodesics on compact manifolds.
\newblock {\em Calc. Var. Partial Differential Equations}, 43(1-2):83--121,
  2012.

\bibitem{MR3913850}
J\"{u}rgen Jost, Enno Ke\ss~ler, J\"{u}rgen Tolksdorf, Ruijun Wu, and Miaomiao
  Zhu.
\newblock From harmonic maps to the nonlinear supersymmetric sigma model of
  quantum field theory: at the interface of theoretical physics, {R}iemannian
  geometry, and nonlinear analysis.
\newblock {\em Vietnam J. Math.}, 47(1):39--67, 2019.

\bibitem{MR3724759}
J\"urgen Jost, Lei Liu, and Miaomiao Zhu.
\newblock A global weak solution of the {D}irac-harmonic map flow.
\newblock {\em Ann. Inst. H. Poincar\'e Anal. Non Lin\'eaire},
  34(7):1851--1882, 2017.

\bibitem{MR2569270}
J{\"u}rgen Jost, Xiaohuan Mo, and Miaomiao Zhu.
\newblock Some explicit constructions of {D}irac-harmonic maps.
\newblock {\em J. Geom. Phys.}, 59(11):1512--1527, 2009.

\bibitem{MR0241822}
O.~A. Lady{\v{z}}enskaja, V.~A. Solonnikov, and N.~N. Ural{\cprime}ceva.
\newblock {\em Linear and quasilinear equations of parabolic type}.
\newblock Translated from the Russian by S. Smith. Translations of Mathematical
  Monographs, Vol. 23. American Mathematical Society, Providence, R.I., 1967.

\bibitem{MR1031992}
H.~Blaine Lawson, Jr. and Marie-Louise Michelsohn.
\newblock {\em Spin geometry}, volume~38 of {\em Princeton Mathematical
  Series}.
\newblock Princeton University Press, Princeton, NJ, 1989.

\bibitem{MR2431658}
Fanghua Lin and Changyou Wang.
\newblock {\em The analysis of harmonic maps and their heat flows}.
\newblock World Scientific Publishing Co. Pte. Ltd., Hackensack, NJ, 2008.

\bibitem{MR1896863}
Seiki Nishikawa.
\newblock {\em Variational problems in geometry}, volume 205 of {\em
  Translations of Mathematical Monographs}.
\newblock American Mathematical Society, Providence, RI, 2002.
\newblock Translated from the 1998 Japanese original by Kinetsu Abe, Iwanami
  Series in Modern Mathematics.

\bibitem{MR604040}
J.~Sacks and K.~Uhlenbeck.
\newblock The existence of minimal immersions of {$2$}-spheres.
\newblock {\em Ann. of Math. (2)}, 113(1):1--24, 1981.

\bibitem{MR826871}
Michael Struwe.
\newblock On the evolution of harmonic mappings of {R}iemannian surfaces.
\newblock {\em Comment. Math. Helv.}, 60(4):558--581, 1985.

\bibitem{MR965544}
Michael Struwe.
\newblock Heat-flow methods for harmonic maps of surfaces and applications to
  free boundary problems.
\newblock In {\em Partial differential equations ({R}io de {J}aneiro, 1986)},
  volume 1324 of {\em Lecture Notes in Math.}, pages 293--319. Springer,
  Berlin, 1988.

\bibitem{MR2431434}
Michael Struwe.
\newblock {\em Variational methods}, volume~34 of {\em Ergebnisse der
  Mathematik und ihrer Grenzgebiete. 3. Folge. A Series of Modern Surveys in
  Mathematics [Results in Mathematics and Related Areas. 3rd Series. A Series
  of Modern Surveys in Mathematics]}.
\newblock Springer-Verlag, Berlin, fourth edition, 2008.
\newblock Applications to nonlinear partial differential equations and
  Hamiltonian systems.

\bibitem{MR2744149}
Michael~E. Taylor.
\newblock {\em Partial differential equations {III}. {N}onlinear equations},
  volume 117 of {\em Applied Mathematical Sciences}.
\newblock Springer, New York, second edition, 2011.

\bibitem{MR1883901}
Peter Topping.
\newblock Reverse bubbling and nonuniqueness in the harmonic map flow.
\newblock {\em Int. Math. Res. Not.}, 2002.

\bibitem{MR2544729}
Changyou Wang and Deliang Xu.
\newblock Regularity of {D}irac-harmonic maps.
\newblock {\em Int. Math. Res. Not. IMRN}, (20):3759--3792, 2009.

\bibitem{MR3719555}
Johannes Wittmann.
\newblock Short time existence of the heat flow for {D}irac-harmonic maps on
  closed manifolds.
\newblock {\em Calc. Var. Partial Differential Equations}, 56(6):Art. 169, 32,
  2017.

\bibitem{MR2496649}
Ling Yang.
\newblock A structure theorem of {D}irac-harmonic maps between spheres.
\newblock {\em Calc. Var. Partial Differential Equations}, 35(4):409--420,
  2009.

\bibitem{MR2267756}
Liang Zhao.
\newblock Energy identities for {D}irac-harmonic maps.
\newblock {\em Calc. Var. Partial Differential Equations}, 28(1):121--138,
  2007.

\bibitem{MR2506243}
Miaomiao Zhu.
\newblock Regularity for weakly {D}irac-harmonic maps to hypersurfaces.
\newblock {\em Ann. Global Anal. Geom.}, 35(4):405--412, 2009.

\end{thebibliography}

\end{document}